\documentclass[twoside,11pt]{amsart}
\usepackage{amsfonts}

\usepackage{latexsym,amsmath,amsfonts}

\usepackage{amsmath}
\usepackage{latexsym}
\usepackage{amsthm}
\usepackage{amssymb}
\usepackage{graphicx}

\newtheorem{theorem}{Theorem}[section]
\newtheorem{remark}{Remark}[section]
\newtheorem{coro}{Corollary}[section]
\newtheorem{lemma}{Lemma}[section]
\newtheorem{prop}{Proposition}[section]
\numberwithin{equation}{section}

\numberwithin{equation}{section}

\begin{document}

\title{\bf The asymptotic behavior of globally smooth solutions of  bipolar non-isentropic compressible
 Euler-Maxwell system for plasma }

\author{ Shu Wang$^{1}$, Yuehong Feng$^{1}$ and Xin Li$^2$ }

\date{}

\maketitle \markboth{S. Wang, Y. H. Feng and X. Li} {  Bipolar
non-isentropic compressible
 Euler-Maxwell system }

\vspace{-3mm}

\begin{center}
{\small $^1$College of Applied Sciences, Beijing University of Technology, Beijing 100124, China \\[2mm]
$^2$Department of Mathematics and computer Science, Xinyang
Vocational and Technical College, Xinyang, 464000, China \\[2mm]
Email~:wangshu@bjut.edu.cn,\hspace{1mm}fengyuehong001@163.com,\hspace{1mm}lixin91600@163.com}
\end{center}


\vspace{1cm}

\begin{center}
\begin{minipage}{14cm}
    {\bf Abstract.} {\small The bipolar non-isentropic compressible Euler-Maxwell system is
investigated in $R^3$ in the present paper, and the $L^q$ time decay
rate for the global smooth solution is established. It is shown that
the total densities, total temperatures and magnetic field of two
carriers converge to the equilibrium states at the same rate
$(1+t)^{-\frac{3}{2}+\frac{3}{2q}}$ in $L^q$ norm.  But, both the
difference of densities and the difference of temperatures of two
carriers decay at the rate $(1+t)^{-2-\frac{1}{q}}$, and the
velocity and electric field decay at the rate
$(1+t)^{-\frac{3}{2}+\frac{1}{2q}}$. This phenomenon on the charge
transport shows the essential difference between the non-isentropic
unipolar Euler-Maxwell and the bipolar isentropic Euler-Maxwell
system. }
\end{minipage}
\end{center}

\vspace{7mm}


\vspace{5mm}

\noindent {\bf AMS Subject Classification (2000)~:} 35A01, 35L45,
35L60, 35Q35


\vspace{3mm}


\section{Introduction and main results}

The Euler-Maxwell system is used to model and simulate the transport
of charged particles in
plasma\cite{Chen84,Dink05,Jero03,Jero05,YW11}. Usually, it takes the
form of compressible non-isentropic Euler equations forced by the
electromagnetic field, which is governed by the self-consistent
Maxwell equation.  In present paper, we consider the Cauchy problem
for the bipolar non-isentropic Euler-Maxwell system
\begin{equation}
\label{0}
\left\{\begin{aligned}
&\partial_t {n_e} + \nabla\cdot(n_e  u_e )=0, \\
&\partial_t (n_e u_e) +\nabla\cdot (n_e u_e\otimes u_e) +\nabla p_e
= -n_e(E+   u_e \times B)-n_e u_e
, \\
&\partial_t(n_e \mathcal {E}_e)+\nabla\cdot (n_e u_e\mathcal
{E}_e+u_e p_e)=-n_e u_e E-n_e |u_e|^2-n_e(\theta_e-1),
\\
&\partial_t {n_i} + \nabla\cdot(n_i  u_i )=0, \\
&\partial_t (n_i u_i) +\nabla\cdot (n_i u_i\otimes u_i) +\nabla p_i
= n_i(E+   u_i \times B)-n_i u_i
, \\
&\partial_t(n_i \mathcal {E}_i)+\nabla\cdot (n_i u_i\mathcal
{E}_i+u_i p_i)=n_i u_i E-n_i |u_i|^2-n_i(\theta_i-1),
\\
&  \partial_t E-\nabla\times B=   n_e u_e-n_i u_i ,
 \\
&  \partial_t B+\nabla\times E=0,\\
&     \nabla\cdot E=n_i-n_e ,\quad \nabla\cdot B=0,\quad
(t,x)\in(0,\infty)\times\mathbb{R}^3,
\end{aligned} \right.
\end{equation}
where the unknowns are the density $n_\mu>0$, the velocity
$u_\mu=(u_\mu^1,u_\mu^2,u_\mu^3)$, the  absolute temperature
$\theta_\mu>0$, the total energy $\mathcal
{E}_\mu=\displaystyle\frac{1}{2}|u_\mu|^2+C_\nu\theta_\mu$, the
pressure function $p_\mu=R_\mu n_\mu\theta_\mu$ for $\mu=e,i$, the
electronic field $E$ and magnetic field $B$. Furthermore, the
constants $C_\nu>0$, $R_\nu>0$  are the heat capacity at constant
volume and the coefficient of heat conductivity respectively.
Throughout this paper, we set $C_\nu=R_\nu=1$ without loss of
generality. Then, the system \eqref{0} is equivalent to
\begin{equation}
\label{1.1}
\left\{\begin{aligned}
&\partial_t {n_e} + \nabla\cdot(n_e  u_e )=0, \\
&\partial_t u_e +(u_e\cdot\nabla) u_e +\frac{\theta_e}{n_e}\nabla
n_e+\nabla\theta_e= - (E+ u_e \times B)-  u_e
, \\
&\partial_t {\theta_e} + \nabla\cdot(\theta_e  u_e )+(\theta_e-1)=0,
\\
&\partial_t {n_i} + \nabla\cdot(n_i  u_i )=0, \\
&\partial_t u_i +(u_i\cdot\nabla) u_i +\frac{\theta_i}{n_i}\nabla
n_i+\nabla\theta_i=  (E+ u_i \times B)-  u_i
, \\
&\partial_t {\theta_i} + \nabla\cdot(\theta_i  u_i )+(\theta_i-1)=0,
\\
&  \partial_t E-\nabla\times B=   n_e u_e-n_i u_i ,
 \\
&  \partial_t B+\nabla\times E=0,\\
&     \nabla\cdot E=n_i-n_e ,\quad \nabla\cdot B=0,\quad
(t,x)\in(0,\infty)\times\mathbb{R}^3.
\end{aligned} \right.
\end{equation}
Initial data is given as
\begin{equation}
\label{1.2}
    (n_\mu, u_\mu, \theta_\mu, E,B)|_{t=0} = (n_{\mu0}, u_{\mu0},
     \theta_{\mu0}, E_0,B_0),\quad x\in
    \mathbb{R}^3,
\end{equation}
with the compatible condition
\begin{equation}
   \label{1.3}
    \nabla\cdot E_0 = n_{i0} - n_{e0}, \quad \nabla\cdot B_0 = 0, \quad x \in \;
    \mathbb{R}^3.
\end{equation}

The Euler-Maxwell system (\ref{1.1}) is a symmetrizable hyperbolic
system for $n_\mu,\theta_\mu >0$. Then the Cauchy problem
(\ref{1.1})-(\ref{1.2}) has a local smooth solution when the initial
data are smooth. In a simplified one dimensional isentropic
Euler-Maxwell system, the global existence of entropy solutions has
been given in \cite{{CJW00}} by the compensated compactness method.
For the three dimensional isentropic Euler-Maxwell system, the
existence of global smooth solutions with small amplitude to the
Cauchy problem in the whole space and to the periodic problem in the
torus is established by Peng et al in \cite{PWG11} and Ueda et al in
\cite{UKW10} respectively, and the decay rate of the smooth solution
when t goes to infinity is obtained by Duan in \cite{Duan11} and
Ueda et al in \cite{UK11}. For asymptotic limits with small
parameters, see \cite{PW08a,PW08b} and references therein. For the
three dimensional bipolar isentropic Euler-Maxwell system, the
global existence and the asymptotic behavior of the smooth solution
is also obtained by Duan et al in \cite{Duan11b}. Recently, Yang et
al in \cite{YW11} consider the diffusive relaxation limit of the
three dimensional unipolar non-isentropic Euler-Maxwell system, and
Wang et al asymptotics and global existence in \cite{FWK11}.

However, there is no analysis on the asymptotics and global
existence for the bipolar non-isentropic Euler-Maxwell system in
three space dimensions yet. Therefore, the goal of the present paper
is to establish the global existence of smooth solutions around a
equilibrium solution of system (\ref{1.1}) and the decay rate of the
smooth solution as $t\rightarrow\infty$.

The main result of this paper can be stated as follows.
\begin{theorem}\label{thm1.1}
Assume \eqref{1.3} hold. If $\left\|
{[n_{\mu0}-1,u_{\mu0},\theta_{\mu0}-1,E_0,B_0] } \right\|_s\leq
\delta_0$ for $s\geq 4$. Then, there is a unique global solution
$[n_\mu(t,x)$, $u_\mu (t,x),$ $\theta_\mu(t,x), E(t,x), B(t,x)]$  to
the initial value problem \eqref{1.1}- \eqref{1.2} which satisfies
 $$[n_\mu -1,u_\mu,\theta_\mu-1, E, B]
 \in C^1\big([0,T);H^{s-1}(\mathbb{R}^3)\big) \cap
C\big([0,T);H^s(\mathbb{R}^3)\big)$$ and
 $$
\mathop {\sup }\limits_{t \geqslant 0}\left\|
{[n_\mu(t)-1,u_\mu(t),\theta_\mu(t)-1,E(t),B(t)] } \right\|_s\leq
C_0\left\| {[n_{\mu0}-1,u_{\mu0},\theta_{\mu0}-1,E_0,B_0] }
\right\|_s,
$$
where $\delta_0, C_0 >0$  are constants independent of time.

Moreover, if $\left\| {[n_{\mu0}-1,u_{\mu0},\theta_{\mu0}-1,E_0,B_0]
} \right\|_{L^1 \cap H^{13} } \leq \delta_1$,
 then the solution
$[n_\mu (t,x)$, $u_\mu (t,x)$,
 $\theta_\mu(t,x)$, $E(t,x)$, $B(t,x)]$ satisfies
\begin{equation}\label{1.4}
 \left\|\left[ {n_e(t)-n_i(t), \theta_e(t)-\theta_i(t) }\right]
\right\|_{L^q}\leq C_1 (1+t)^{-2-\frac{1}{q}},
\end{equation}
\begin{equation}\label{1.5}
\left\|\left[ {n_e(t)+n_i(t)-2, \theta_e(t)+\theta_i(t)-2 }\right]
\right\|_{L^q}\leq C_1 (1+t)^{-\frac{3}{2}+\frac{3}{2q}},
\end{equation}
\begin{equation}\label{1.6}
\left\| {u_e(t)\pm u_i(t),  E(t)  } \right\|_{L^q}\leq C_1
(1+t)^{-\frac{3}{2}+\frac{1}{2q}},
 \end{equation}
\begin{equation}\label{1.7}
  \left\|
{B}(t)  \right\|_{L^q}\leq C_1 (1+t)^{-\frac{3}{2}+\frac{3}{2q}},
\end{equation} for any $t\geq 0$ and $2\leq q \leq\infty$.
Where, constants $\delta_1, C_1 >0$ are also independent of time.
\end{theorem}
\begin{remark}It should be emphasized that both the velocity and temperature relaxation term of the
bipolar non-isentropic Euler-Maxwell system \eqref{1.1} plays a key
role in the proof of Theorem \ref{thm1.1}.
\end{remark}

\noindent {\bf Notations.}
 In this paper, $f \sim
g$ means $\gamma a\leq b\leq \frac{1}{\gamma}$ for a constant
$0<\gamma<1$. $H^s$ denotes the standard Sobolev space
$W^{s,2}({{\mathbb R}^3})$. We use $\dot{H}^s$ to denote the
 corresponding $s$-order homogeneous Sobolev space. Set $L^2=H^0$. The
norm of $H^s$ is denoted by $\left\| {\cdot} \right\|_s$ with
$\left\| {\cdot} \right\|=\left\| {\cdot} \right\|_0$, and
$\langle\cdot,\cdot\rangle$  denotes the inner product over
$L^2({\mathbb R}^3)$.  For the multi-index $\alpha
=(\alpha_1,\alpha_2,\alpha_3)$, we denote $\partial^\alpha$
$=\partial_{x_1}^{\alpha_1}
\partial_{x_2}^{\alpha_2}\partial_{x_3}^{\alpha_3}$ $=\partial_{ 1}^{\alpha_1}
\partial_{ 2}^{\alpha_2}\partial_{ 3}^{\alpha_3}$ and
$|\alpha|=\alpha_1+\alpha_2+\alpha_3$. For an integrable function
$f:{\mathbb R}^3\rightarrow{\mathbb R}$, its Fourier transform is
defined by
$$ \hat{f}(k)=\int_{{\mathbb R}^3}e^{-ix\cdot k}f(x)dx,\ \ x\cdot k:=
\sum\limits_{j=1}^{3} {x_j k_j},\ \ k\in {\mathbb R}^3,
$$
where $i=\sqrt{-1}\in \mathbb C$ is the imaginary unit.

The rest of the paper is arranged as follows. In Section 2, the
transformation of the initial value problem and the proof of the
global existence and uniqueness of solutions are presented. In
Section 3, we study the linearized homogeneous equations to get the
$L^p-L^q$ decay property and the explicit representation of
solutions. In the last Section 4, we investigate the decay rates of
solutions to the transformed nonlinear equations and complete the
proof of Theorem \ref{thm1.1}.

\section{Global solutions for equations \eqref{1.1}}
\subsection{Preliminary}
Suppose $[n_\mu (t,x), u_\mu (t,x),\theta_\mu (t,x),E(t,x),B(t,x)]$
be a smooth solution of the initial value problem for the bipolar
non-isentropic Euler-Maxwell equations \eqref{1.1} with initial data
\eqref{1.2} which satisfies \eqref{1.3}. Set
\begin{equation}\label{2.1}
n_\mu (t,x)=1+\rho_\mu (t,x), \theta_\mu (t,x)=1+\Theta_\mu (t,x).
\end{equation}
Thus, we can rewrite the system \eqref{1.1}-\eqref{1.3} as
\begin{equation}
\label{2.2}
\left\{\begin{aligned}
&\partial_t {\rho_e} + \nabla\cdot((1+\rho_e)  u_e )=0, \\
&\partial_t u_e +(u_e\cdot\nabla) u_e
+\frac{1+\Theta_e}{1+\rho_e}\nabla \rho_e +\nabla\Theta_e= - (E+ u_e
\times B)- u_e
, \\
&\partial_t {\Theta_e} + \nabla\cdot((1+\Theta_e)  u_e )+\Theta_e=0,
\\
&\partial_t {\rho_i} + \nabla\cdot((1+\rho_i)  u_i )=0, \\
&\partial_t u_i +(u_i\cdot\nabla) u_i
+\frac{1+\Theta_i}{1+\rho_i}\nabla \rho_i+\nabla\Theta_i=  (E+ u_i
\times B)-  u_i
, \\
&\partial_t {\Theta_i} + \nabla\cdot((1+\Theta_i)  u_i )+\Theta_i=0,
\\
&  \partial_t E-\nabla\times B -u_e+ u_i=   \rho_e u_e-\rho_i u_i ,
 \\
&  \partial_t B+\nabla\times E=0,\\
&     \nabla\cdot E=\rho_i-\rho_e ,\quad \nabla\cdot B=0,\quad
(t,x)\in(0,\infty)\times\mathbb{R}^3,
\end{aligned} \right.
\end{equation}
with initial data
\begin{equation}\label{2.3}
U|_{t=0}=U_0:=[\rho_{\mu 0},u_{\mu 0},\Theta_{ \mu0},E_0,B_0],\ x
\in{\mathbb R}^3,
\end{equation}
which satisfies the compatible condition
\begin{equation}\label{2.4}
    \nabla\cdot E_0 = \rho_{i0}-\rho_{e0}, \quad \nabla\cdot B_0 = 0, \quad x \in  \;
    \mathbb{R}^3.
\end{equation}
Here, $\rho_{\mu 0}=n_{\mu 0}-1.$

In the following, we usually assume $s\geq4$. Moreover, for
$U=[\rho_{\mu }$, $u_{\mu }$, $\Theta_{\mu }$, $E$, $B]$, we use
$\mathcal {E}_s(U(t))$, $\mathcal {E}_s^h(U(t))$, $\mathcal
{D}_s(U(t))$ and $\mathcal {D}_s^h(U(t))$ to define the energy
functional, the high-order energy functional, the dissipation rate
and the high-order dissipation rate as
\begin{equation}\label{2.5}
\mathcal {E}_s(U(t))\sim \left\| {[\rho_\mu,u_\mu,\Theta_\mu,E,B]}
\right\|_s^2,
\end{equation}
\begin{equation}\label{2.6}
\mathcal {E}_s^h(U(t))\sim \left\|
{\nabla[\rho_\mu,u_\mu,\Theta_\mu,E,B]} \right\|_{s-1}^2,
\end{equation}
\begin{equation}\label{2.7}
\begin{split}
\mathcal {D}_s(U(t))\sim & \left\| {\nabla[\rho_e,\rho_i]}
\right\|_{s-1}^2 +\left\| {[u_e,u_i,\Theta_e,\Theta_i]}
\right\|_{s}^2\\
&+\left\| {E} \right\|_{s-1}^2+\left\| {\nabla B}
\right\|_{s-2}^2+\left\| {\rho_e-\rho_i} \right\|^2
\end{split}
\end{equation}
and
\begin{equation}\label{2.8}
\begin{split}
\mathcal {D}^h_s(U(t))\sim & \left\| {\nabla^2[\rho_e,\rho_i]}
\right\|_{s-2}^2 +\left\| {\nabla[u_e,u_i,\Theta_e,\Theta_i]}
\right\|_{s-1}^2\\
&+\left\| {\nabla E} \right\|_{s-2}^2+\left\| {\nabla^2 B}
\right\|_{s-3}^2+\left\| {\nabla [\rho_e-\rho_i]} \right\|^2,
\end{split}\end{equation}
 respectively. Now, concerning the transformed initial value problem
\eqref{2.2}-\eqref{2.3}, we have the global existence result as
follows.
\begin{prop}\label{prop2.1}
Assume that $U_0=[\rho_{\mu0},u_{\mu0},\Theta_{\mu0},E_0,B_0]$
satisfies the compatible condition \eqref{2.4}. If \ $ \mathcal
{E}_s(U_0) $ is small enough, then, for any $t \geq0$, the initial
value problem \eqref{2.2}-\eqref{2.3} has a unique global nonzero
solution $U =[\rho_\mu,u_\mu,\Theta_\mu,E,B]$ which satisfies
\begin{equation}\label{2.9}
U\in C^1\big([0,T);H^{s-1}(\mathbb{R}^3)\big) \cap
C\big([0,T);H^s(\mathbb{R}^3)\big),
\end{equation}
and
\begin{equation}\label{2.10}
\mathcal {E}_s(U(t))+\lambda\int_0^t\mathcal {D}_s(U(s))ds\leq
\mathcal {E}_s(U_0).
\end{equation}
\end{prop}
Obviously, from the Proposition \ref{prop2.1}, it is straightforward
to get the existence result of Theorem \ref{thm1.1}. Furthermore,
solutions of Proposition \ref{prop2.1} really decay under some extra
conditions on $U_0=[\rho_{\mu0},u_{\mu0},\Theta_{\mu0},E_0,B_0]$.
For this purpose, we define $\omega_s(U_0)$ as
\begin{equation}\label{2.11}
\omega_s(U_0)=\left\|{U_0}\right\|_s+\left\|{[\rho_{\mu0}, u_{\mu0},
\Theta_{\mu0}, E_0,B_0]}\right\|_{L^1}
\end{equation}
for $s\geq4.$ Then, we obtain the following decay results.
\begin{prop}\label{prop2.2}
Assume that $U_0=[\rho_{\mu0}$, $u_{\mu0}$, $\Theta_{\mu0}$, $E_0$,
$B_0]$ satisfies \eqref{2.4}. If $ \omega_{s+2}(U_0)$ is
sufficiently small, then system \eqref{2.2}-\eqref{2.4} has a
solution $U =[\rho_\mu$, $u_\mu$, $\Theta_\mu$, $E$, $B ]$
satisfying
\begin{equation}\label{2.12}
\left\|{U(t)}\right\|_s\leq C \omega_{s+2}(U_0)(1+t)^{-\frac{3}{4}}
\end{equation}
for any $t\geq0$. Moreover, if $\omega_{s+6}(U_0)$ is sufficiently
small, then, for any $t \geq0$, the solution also satisfies
\begin{equation}\label{2.13}
\left\|{\nabla U(t)}\right\|_{s-1}\leq C
 \omega_{s+6}(U_0)(1+t)^{-\frac{5}{4}}.
\end{equation}
\end{prop}
Thus, one can obtain the decay rates  \eqref{1.4}-\eqref{1.7}
through the method of bootstrap and the Proposition stated above.


\subsection{Weighted energy estimates.} In this subsection, we shall give the
proof of Proposition \ref{prop2.1} for the global existence and
uniqueness of solutions to the initial value problem
\eqref{2.2}-\eqref{2.3}. Since hyperbolic equations \eqref{2.2} is
quasi-linear symmetrizable, thus one has the local existence of
smooth solutions to  \eqref{2.2} as follows.
\begin{lemma}
\label{L2.1}
    (Local existence of smooth solutions, see \cite{Ka75,Ma84})
Let $s > \frac{5}{2}$ and $(\rho_{\mu 0}$, $u_{\mu 0}$, $\Theta_{\mu
0}$, $E_0$, $B_0)$ $\in H^s(\mathbb{R}^3)$. Then there exist $T> 0$ and a unique smooth solution
$(n_\mu$, $u_\mu$, $\theta_\mu$, $E$, $B)$  to the Cauchy problem
(\ref{1.1})-(\ref{1.2}) satisfying $(\rho_{\mu}$, $u_{\mu}$, $\Theta_{\mu}$, $E$, $B)\in
C^1\big([0,T);H^{s-1}(\mathbb{R}^3)\big) \cap
C\big([0,T);H^s(\mathbb{R}^3)\big)$.
\end{lemma}
Then, with the help of the continuity argument, the global existence
of solutions satisfying \eqref{2.9} and \eqref{2.10} follows by
combing Lemma \ref{L2.1} and a priori estimate as follows.
\begin{theorem}\label{thm2.1}
  Assume that
$U$ $=[\rho_\mu$, $u_\mu$, $\Theta_\mu$, $E$, $B ]$ $\in
C^1\big([0,T);H^{s-1}(\mathbb{R}^3)\big) \cap
C\big([0,T);H^s(\mathbb{R}^3)\big)$ is smooth for $T>0$ with
\begin{equation}\label{2.14}
\mathop {\sup }\limits_{0\leq t\leq T
}\left\|{U(t)}\right\|_s\leq\delta
\end{equation}
for $\delta\leq \delta_0$ with $\delta_0$ sufficiently small and
suppose $U$ to be the solution of the equations \eqref{2.2} for
$t\in(0,T)$. Then, for a constant $0<\gamma<1$ and any $0\leq t\leq
T$, it holds that
\begin{equation}\label{2.15}
\frac{d}{dt}\mathcal {E}_s(U(t))+\gamma \mathcal {D}_s(U(t))\leq
C[\mathcal {E}_s(U(t))^\frac{1}{2}+\mathcal {E}_s(U(t))]\mathcal
{D}_s(U(t)).
\end{equation}
\end{theorem}
\noindent \emph{Proof.} We will use five steps to finish the proof
as follows. In step 1, we establish the estimate of Euler part and
Maxwell part of the system \eqref{2.2} by using weighted energy
estimate method. In the following steps $2-4$, we utilize the
skew-symmetric structure of the system \eqref{2.2} to get the
dissipative estimates for $\rho_\mu$, $E$ and $B$. \\
\emph{Step 1.}
It holds that
\begin{equation}\label{2.16}
\frac{d} {{dt}}\left\| U \right\|_s^2 + \left\| {\left[
{{u_e},{u_i},{\Theta _e},{\Theta _i}} \right]} \right\|_s^2
\leqslant C{\left\| U \right\|_s}\left( {\left\| {\left[
{{u_e},{u_i},{\Theta _e},{\Theta _i}} \right]} \right\|_s^2 +
\left\| {\nabla \left[ {{\rho _e},{\rho _i}} \right]} \right\|_{s -
1}^2} \right).
\end{equation}
In fact, from the first six equations of \eqref{2.2}, weighted
energy estimate on $\partial^\alpha \rho_\mu$, $\partial^\alpha
u_\mu$ and $\partial^\alpha \Theta_\mu$ with $|\alpha|\leq s$ imply
\begin{equation}\label{2.17}
\begin{split}  \frac{1} {2}&\frac{d} {{dt}} \sum\limits_{\mu  = e,i} {\left(
{\left\langle {\frac{{1 + {\Theta _\mu }}} {{1 + {\rho _\mu
}}},{{\left| {{\partial ^\alpha }{\rho _\mu }} \right|}^2}}
\right\rangle  + \left\langle {1 + {\rho _\mu },{{\left| {{\partial
^\alpha }{u_\mu }} \right|}^2}} \right\rangle  + \left\langle
{\frac{{1 + {\rho _\mu }}}
{{1 + {\Theta _\mu }}},{{\left| {{\partial ^\alpha }{\Theta _\mu }} \right|}^2}} \right\rangle } \right)}   \\
 &  + \sum\limits_{\mu  = e,i} {\left( {\left\langle {1 + {\rho _\mu },{{\left| {{\partial ^\alpha }{u_\mu }}
   \right|}^2}} \right\rangle  + \left\langle {\frac{{1 + {\rho _\mu }}}
{{1 + {\Theta _\mu }}},{{\left| {{\partial ^\alpha }{\Theta _\mu }} \right|}^2}} \right\rangle } \right)}
+ \left\langle {\left( {1 + {\rho _e}} \right){\partial ^\alpha }E,{\partial ^\alpha }{u_e}} \right\rangle   \\
&   - \left\langle {\left( {1 + {\rho _i}} \right){\partial ^\alpha
}E,{\partial ^\alpha }{u_i}} \right\rangle
   =  - \sum\limits_{\beta  < \alpha } {C_\beta ^\alpha {I_{\alpha ,\beta }}(t) + {I_1}(t).}
\end{split}\end{equation}
Where, ${I_{\alpha ,\beta }}(t)={I^e_{\alpha ,\beta
}}(t)+{I^i_{\alpha ,\beta }}(t)$, $I_1(t)=I_1^e(t)+I_1^i(t)$ with
\begin{equation}\notag\begin{split}
 I_{\alpha ,\beta }^e(t) =&  \left\langle {\frac{{1 + {\Theta _e}}}
{{1 + {\rho _e}}}{\partial ^{\alpha  - \beta }}{\rho _e}\nabla
{\partial ^\beta }{u_e},{\partial ^\alpha }{\rho _e}} \right\rangle
+ \left\langle {\frac{{1 + {\Theta _e}}} {{1 + {\rho _e}}}{\partial
^{\alpha  - \beta }}{u_e}\nabla {\partial ^\beta }
{\rho _e},{\partial ^\alpha }{\rho _e}} \right\rangle   \\
&   + \left\langle {\frac{{1 + {\rho _e}}} {{1 + {\Theta
_e}}}{\partial ^{\alpha  - \beta }}{u_e}\nabla {\partial ^\beta
}{\Theta _e},{\partial ^\alpha }{\Theta _e}} \right\rangle  +
\left\langle {\frac{{1 + {\rho _e}}} {{1 + {\Theta _e}}}{\partial
^{\alpha  - \beta }}{\Theta _e}\nabla {\partial ^\beta }{u_e},
{\partial ^\alpha }{\Theta _e}} \right\rangle   \\
 &  + \left\langle {\left( {1 + {\rho _e}} \right){\partial ^{\alpha  - \beta }}{u_e}\nabla
    {\partial ^\beta }{u_e},{\partial ^\alpha }{u_e}} \right\rangle  + \left\langle
    {\left( {1 + {\rho _e}} \right){\partial ^{\alpha  - \beta }}\left( {\frac{{1 + {\Theta _e}}}
{{1 + {\rho _e}}}} \right)\nabla {\partial ^\beta }{\rho _e},{\partial ^\alpha }{u_e}} \right\rangle   \\
  & + \left\langle {\left( {1 + {\rho _e}} \right){\partial ^{\alpha  - \beta }}{u_e} \times
   {\partial ^\beta }B,{\partial ^\alpha }{u_e}} \right\rangle,
\end{split}\end{equation}
\begin{equation}\notag\begin{split}
 I_{\alpha ,\beta }^i(t) =&  \left\langle {\frac{{1 + {\Theta _i}}}
{{1 + {\rho _i}}}{\partial ^{\alpha  - \beta }}{\rho _i}\nabla
{\partial ^\beta }{u_i},{\partial ^\alpha }{\rho _i}} \right\rangle
+ \left\langle {\frac{{1 + {\Theta _i}}} {{1 + {\rho _i}}}{\partial
^{\alpha  - \beta }}{u_i}\nabla {\partial ^\beta }
{\rho _i},{\partial ^\alpha }{\rho _i}} \right\rangle   \\
&   + \left\langle {\frac{{1 + {\rho _i}}} {{1 + {\Theta
_i}}}{\partial ^{\alpha  - \beta }}{u_i}\nabla {\partial ^\beta
}{\Theta _i},{\partial ^\alpha }{\Theta _i}} \right\rangle  +
\left\langle {\frac{{1 + {\rho _i}}} {{1 + {\Theta _i}}}{\partial
^{\alpha  - \beta }}{\Theta _i}\nabla {\partial ^\beta }{u_i},
{\partial ^\alpha }{\Theta _i}} \right\rangle   \\
 &  + \left\langle {\left( {1 + {\rho _i}} \right){\partial ^{\alpha  - \beta }}{u_i}\nabla
    {\partial ^\beta }{u_i},{\partial ^\alpha }{u_i}} \right\rangle  + \left\langle
    {\left( {1 + {\rho _i}} \right){\partial ^{\alpha  - \beta }}\left( {\frac{{1 + {\Theta _i}}}
{{1 + {\rho _i}}}} \right)\nabla {\partial ^\beta }{\rho _i},{\partial ^\alpha }{u_i}} \right\rangle   \\
  & - \left\langle {\left( {1 + {\rho _i}} \right){\partial ^{\alpha  - \beta }}{u_i} \times
   {\partial ^\beta }B,{\partial ^\alpha }{u_i}} \right\rangle,
\end{split}\end{equation}
and
\begin{equation}\notag\begin{split}
  I_1^e(t) = &~~\frac{1}
{2}\left\langle {{\partial _t}\left( {\frac{{1 + {\Theta _e}}}
{{1 + {\rho _e}}}} \right),{{\left| {{\partial ^\alpha }{\rho _e}} \right|}^2}}
 \right\rangle  + \left\langle {\nabla {\Theta _e}{\partial ^\alpha }{u_e},{\partial ^\alpha }{\rho _e}}
  \right\rangle  + \left\langle {\nabla {\rho _e}{\partial ^\alpha }{u_e},{\partial ^\alpha }{\Theta _e}}
  \right\rangle  \\
  & + \frac{1}
{2}\left\langle {\nabla  \cdot \left( {\frac{{1 + {\Theta _e}}} {{1
+ {\rho _e}}}{u_e}} \right),{{\left| {{\partial ^\alpha }{\rho _e}}
\right|}^2}} \right\rangle  + \frac{1} {2}\left\langle {{\partial
_t}\left( {\frac{{1 + {\rho _e}}}
{{1 + {\Theta _e}}}} \right),{{\left| {{\partial ^\alpha }{\Theta _e}} \right|}^2}} \right\rangle    \\
  & + \frac{1}
{2}\left\langle {\nabla  \cdot \left( {\frac{{1 + {\rho _e}}}
{{1 + {\Theta _e}}}{u_e}} \right),{{\left| {{\partial ^\alpha }{\Theta _e}} \right|}^2}} \right\rangle
 - \left\langle {\left( {1 + {\rho _e}} \right){u_e} \times {\partial ^\alpha }B,{\partial ^\alpha }{u_e}}
 \right\rangle ,
\end{split}\end{equation}
\begin{equation}\notag\begin{split}
  I_1^i(t) = &~~\frac{1}
{2}\left\langle {{\partial _t}\left( {\frac{{1 + {\Theta _i}}} {{1 +
{\rho _i}}}} \right),{{\left| {{\partial ^\alpha }{\rho _i}}
\right|}^2}}
 \right\rangle  + \left\langle {\nabla {\Theta _i}{\partial ^\alpha }{u_i},{\partial ^\alpha }{\rho _i}}
  \right\rangle  + \left\langle {\nabla {\rho _i}{\partial ^\alpha }{u_i},{\partial ^\alpha }{\Theta _i}}
  \right\rangle  \\
  & + \frac{1}
{2}\left\langle {\nabla  \cdot \left( {\frac{{1 + {\Theta _i}}} {{1
+ {\rho _i}}}{u_i}} \right),{{\left| {{\partial ^\alpha }{\rho _i}}
\right|}^2}} \right\rangle  + \frac{1} {2}\left\langle {{\partial
_t}\left( {\frac{{1 + {\rho _i}}}
{{1 + {\Theta _i}}}} \right),{{\left| {{\partial ^\alpha }{\Theta _i}} \right|}^2}} \right\rangle    \\
  & + \frac{1}
{2}\left\langle {\nabla  \cdot \left( {\frac{{1 + {\rho _i}}} {{1 +
{\Theta _i}}}{u_i}} \right),{{\left| {{\partial ^\alpha }{\Theta
_i}} \right|}^2}} \right\rangle
 +\left\langle {\left( {1 + {\rho _i}} \right){u_i} \times {\partial ^\alpha }B,{\partial ^\alpha }{u_i}}
 \right\rangle ,
\end{split}\end{equation}
where we have used integration by parts. When $|\alpha|=0$, one has
\begin{equation}\notag\begin{split}
  {I_1}(t) =& I_1^e(t) + I_1^i(t)  \\
   = &\sum\limits_{\mu  = e,i} {\left( {\frac{1}
{2}\left\langle {{\partial _{{\Theta _\mu }}}\left( {\frac{{1 +
{\Theta _\mu }}} {{1 + {\rho _\mu }}}} \right){\partial _t}{\Theta
_\mu } + {\partial _{{\rho _\mu }}}\left( {\frac{{1 + {\Theta _\mu
}}} {{1 + {\rho _\mu }}}} \right){\partial _t}{\rho _\mu }, {{\left|
{{\rho _\mu }} \right|}^2}} \right\rangle
+ \left\langle {\nabla {\Theta _\mu }{u_\mu },{\rho _\mu }} \right\rangle } \right.}   \\
 &  + \frac{1}
{2}\left\langle {{\partial _{{\Theta _\mu }}}\left( {\frac{{1 +
{\Theta _\mu }}} {{1 + {\rho _\mu }}}{u_\mu }} \right)\nabla {\Theta
_\mu } + {\partial _{{u_\mu }}}\left( {\frac{{1 + {\Theta _\mu }}}
{{1 + {\rho _\mu }}}{u_\mu }} \right)\nabla  \cdot {u_\mu } +
{\partial _{{\rho _\mu }}}\left( {\frac{{1 + {\Theta _\mu }}}
{{1 + {\rho _\mu }}}{u_\mu }} \right)\nabla {\rho _\mu },{{\left| {{\rho _\mu }} \right|}^2}} \right\rangle   \\
 &  + \frac{1}
{2}\left\langle {{\partial _{{\Theta _\mu }}}\left( {\frac{{1 +
{\rho _\mu }}} {{1 + {\Theta _\mu }}}} \right){\partial _t}{\Theta
_\mu } + {\partial _{{\rho _\mu }}}\left( {\frac{{1 + {\rho _\mu }}}
{{1 + {\Theta _\mu }}}} \right){\partial _t}{\rho _\mu },{{\left| {{\Theta _\mu }} \right|}^2}} \right\rangle
  + \left\langle {\nabla {\rho _\mu }{ }{u_\mu },{ }{\Theta _\mu }} \right\rangle  \\
 & \left. { + \frac{1}
{2}\left\langle {{\partial _{{\Theta _\mu }}}\left( {\frac{{1 +
{\rho _\mu }}} {{1 + {\Theta _\mu }}}{u_\mu }} \right)\nabla {\Theta
_\mu } + {\partial _{{u_\mu }}}\left( {\frac{{1 + {\rho _\mu }}} {{1
+ {\Theta _\mu }}}{u_\mu }} \right)\nabla  \cdot {u_\mu } +
{\partial _{{\rho _\mu }}}\left( {\frac{{1 + {\rho _\mu }}}
{{1 + {\Theta _\mu }}}{u_\mu }} \right)\nabla {\rho _\mu },{{\left| {{\Theta _\mu }}
 \right|}^2}} \right\rangle } \right) \\
   &- \left\langle {\left( {1 + {\rho _e}} \right){u_e} \times B,{u_e}} \right\rangle
   +\left\langle {\left( {1 + {\rho _i}} \right){u_i} \times B,{u_i}}
   \right\rangle
\end{split}\end{equation}
   \begin{equation}\notag\begin{split}
   =& \sum\limits_{\mu  = e,i} {\left( { - \frac{1}
{2}\left\langle {{\partial _{{\Theta _\mu }}}\left( {\frac{{1 +
{\Theta _\mu }}} {{1 + {\rho _\mu }}}} \right)\nabla  \cdot \left(
{{u_\mu }\left( {1 + {\Theta _\mu }} \right)} \right) + {\partial
_{{\rho _\mu }}}\left( {\frac{{1 + {\Theta _\mu }}} {{1 + {\rho _\mu
}}}} \right)\nabla  \cdot \left( {{u_\mu }\left( {1 + {\rho _\mu }}
\right)} \right),
{{\left| {{\rho _\mu }} \right|}^2}} \right\rangle } \right.}   \\
  & + \frac{1}
{2}\left\langle {{\partial _{{\Theta _\mu }}}\left( {\frac{{1 +
{\Theta _\mu }}} {{1 + {\rho _\mu }}}{u_\mu }} \right)\nabla {\Theta
_\mu } + {\partial _{{u_\mu }}}\left( {\frac{{1 + {\Theta _\mu }}}
{{1 + {\rho _\mu }}}{u_\mu }} \right)\nabla  \cdot {u_\mu } +
{\partial _{{\rho _\mu }}}\left( {\frac{{1 + {\Theta _\mu }}}
{{1 + {\rho _\mu }}}{u_\mu }} \right)\nabla {\rho _\mu },{{\left| {{\rho _\mu }}
\right|}^2}} \right\rangle   \\
  & - \frac{1}
{2}\left\langle {{\partial _{{\Theta _\mu }}}\left( {\frac{{1 +
{\rho _\mu }}} {{1 + {\Theta _\mu }}}} \right)\nabla  \cdot \left(
{{u_\mu }\left( {1 + {\Theta _\mu }} \right)} \right) + {\partial
_{{\rho _\mu }}}\left( {\frac{{1 + {\rho _\mu }}}
{{1 + {\Theta _\mu }}}} \right)\nabla  \cdot \left( {{u_\mu }\left( {1 + {\rho _\mu }}
\right)} \right),{{\left| {{\Theta _\mu }} \right|}^2}} \right\rangle   \\
 &  + \frac{1}
{2}\left\langle {{\partial _{{\Theta _\mu }}}\left( {\frac{{1 +
{\rho _\mu }}} {{1 + {\Theta _\mu }}}{u_\mu }} \right)\nabla {\Theta
_\mu } + {\partial _{{u_\mu }}}\left( {\frac{{1 + {\rho _\mu }}} {{1
+ {\Theta _\mu }}}{u_\mu }} \right)\nabla  \cdot {u_\mu } +
{\partial _{{\rho _\mu }}}\left( {\frac{{1 + {\rho _\mu }}}
{{1 + {\Theta _\mu }}}{u_\mu }} \right)\nabla {\rho _\mu },{{\left| {{\Theta _\mu }} \right|}^2}} \right\rangle   \\
 & \left. { + \left\langle {\nabla {\rho _\mu }{ }{u_\mu },{ }{\Theta _\mu }}
  \right\rangle  + \left\langle {\nabla {\Theta _\mu }{u_\mu },{\rho _\mu }} \right\rangle } \right)
   - \left\langle {\left( {1 + {\rho _e}} \right){u_e} \times B,{u_e}} \right\rangle
   + \left\langle {\left( {1 + {\rho _i}} \right){u_i} \times B,{u_i}} \right\rangle
\end{split}\end{equation}
\begin{equation}\notag\begin{split}
     \leqslant & C\left\| {{\rho _\mu }} \right\|{\left\| {{\rho _\mu }} \right\|_{{L^\infty }}}\left
     \{ {{{\left\| {{\partial _{{\Theta _\mu }}}\left( {\frac{{1 + {\Theta _\mu }}}
{{1 + {\rho _\mu }}}} \right)} \right\|}_{{L^\infty }}}\left( {{{\left\| {1 + {\Theta _\mu }}
\right\|}_{{L^\infty }}}\left\| {\nabla  \cdot {u_\mu }} \right\| + {{\left\| {\nabla {\Theta _\mu }}
\right\|}_{{L^\infty }}}\left\| {{u_\mu }} \right\|} \right)} \right.  \\
  & + {\left\| {{\partial _{{\rho _\mu }}}\left( {\frac{{1 + {\Theta _\mu }}}
{{1 + {\rho _\mu }}}} \right)} \right\|_{{L^\infty }}}\left(
{{{\left\| {1 + {\rho _\mu }} \right\|}_{{L^\infty }}}\left\|
{\nabla  \cdot {u_\mu }} \right\| + {{\left\| {\nabla {\rho _\mu }}
\right\|}_{{L^\infty }}}\left\| {{u_\mu }} \right\|} \right) +
{\left\| {{\partial _{{\Theta _\mu }}}\left( {\frac{{1 + {\Theta
_\mu }}}
{{1 + {\rho _\mu }}}{u_\mu }} \right)} \right\|_{{L^\infty }}}   \\
  &  \left. {\left\| {\nabla {\Theta _\mu }} \right\| + {{\left\| {{\partial _{{u_\mu }}}\left(
  {\frac{{1 + {\Theta _\mu }}}
{{1 + {\rho _\mu }}}{u_\mu }} \right)} \right\|}_{{L^\infty
}}}\left\| {\nabla  \cdot {u_\mu }} \right\| + {{\left\| {{\partial
_{{\rho _\mu }}}\left( {\frac{{1 + {\Theta _\mu }}}
{{1 + {\rho _\mu }}}{u_\mu }} \right)} \right\|}_{{L^\infty }}}\left\| {\nabla {\rho _\mu }} \right\|} \right\}  \\
  & + C\left\| {{\Theta _\mu }} \right\|{\left\| {{\Theta _\mu }} \right\|_{{L^\infty }}}\left\{ {{{\left\|
  {{\partial _{{\Theta _\mu }}}\left( {\frac{{1 + {\rho _\mu }}}
{{1 + {\Theta _\mu }}}} \right)} \right\|}_{{L^\infty }}}\left( {{{\left\| {1 + {\Theta _\mu }}
 \right\|}_{{L^\infty }}}\left\| {\nabla  \cdot {u_\mu }} \right\| + {{\left\| {\nabla {\Theta _\mu }}
 \right\|}_{{L^\infty }}}\left\| {{u_\mu }} \right\|} \right)} \right. \\
 &  + {\left\| {{\partial _{{\rho _\mu }}}\left( {\frac{{1 + {\rho _\mu }}}
{{1 + {\Theta _\mu }}}} \right)} \right\|_{{L^\infty }}}\left(
{{{\left\| {1 + {\rho _\mu }} \right\|}_{{L^\infty }}}\left\|
{\nabla  \cdot {u_\mu }} \right\| + {{\left\| {\nabla {\rho _\mu }}
\right\|}_{{L^\infty }}}\left\| {{u_\mu }} \right\|} \right) +
{\left\| {{\partial _{{\Theta _\mu }}}\left( {\frac{{1 + {\rho _\mu
}}}
{{1 + {\Theta _\mu }}}{u_\mu }} \right)} \right\|_{{L^\infty }}}    \\
 & \left. {  \left\| {\nabla {\Theta _\mu }} \right\| + {{\left\| {{\partial _{{u_\mu }}}\left(
  {\frac{{1 + {\rho _\mu }}}
{{1 + {\Theta _\mu }}}{u_\mu }} \right)} \right\|}_{{L^\infty
}}}\left\| {\nabla  \cdot {u_\mu }} \right\| + {{\left\| {{\partial
_{{\rho _\mu }}}\left( {\frac{{1 + {\rho _\mu }}}
{{1 + {\Theta _\mu }}}{u_\mu }} \right)} \right\|}_{{L^\infty }}}\left\|
{\nabla {\rho _\mu }} \right\|} \right\}  \\
  & + C\left\| {\nabla {\rho _\mu }} \right\|\left\| {{u_\mu }} \right\|{\left\|
   {{\Theta _\mu }} \right\|_{{L^\infty }}} + C\left\| {\nabla {\Theta _\mu }}
   \right\|\left\| {{u_\mu }} \right\|{\left\| {{\rho _\mu }} \right\|_{{L^\infty }}}
   + C{\left\| {1 + {\rho _\mu }} \right\|_{{L^\infty }}}\left\| {{u_\mu }} \right\|\left\| B \right\|{
   \left\| {{u_\mu }} \right\|_{{L^\infty }}}
\end{split}\end{equation}
\begin{equation}\notag\begin{split}
   \leqslant & C\left( {\left\| {\nabla {u_\mu }} \right\| + \left\| {{u_\mu }} \right\| +
    \left\| {\nabla {\Theta _\mu }} \right\| + \left\| {\nabla {\rho _\mu }} \right\|} \right)\left(
    {\left\| {{\rho _\mu }} \right\|{{\left\| {\nabla {\rho _\mu }} \right\|}_1} + \left\| {{\Theta _\mu }}
     \right\|{{\left\| {\nabla {\Theta _\mu }} \right\|}_1}} \right)   \\
   & + C\left\| {\nabla {\rho _\mu }} \right\|\left\| {{u_\mu }} \right\|\left\| {\nabla {\Theta _\mu }} \right\|
    + \left\| {\nabla {\Theta _\mu }} \right\|\left\| {{u_\mu }} \right\|{\left\| {\nabla {\rho _\mu }} \right\|_1}
    + C\left\| {{u_\mu }} \right\|\left\| B \right\|{\left\| {\nabla {u_\mu }} \right\|_1}  \\
   \leqslant & C\left\| {\left[ {{\rho _\mu },{u_\mu },{\Theta _\mu },B} \right]} \right\|\left( {\left\|
    {\nabla {\rho _\mu }} \right\|_1^2 + \left\| {{u_\mu }} \right\|_2^2 + \left\| {\nabla {\Theta _\mu }}
     \right\|_1^2} \right),
\end{split}\end{equation}
which will further be bounded by the right hand side term of
\eqref{2.16}, and where we have used \eqref{2.14}. When
$|\alpha|\geq 1$, similarly as before, one has
\[{I_{\alpha ,\beta }}(t) + {I_1}(t) \leqslant C{\left\| {\left[ {{\rho _\mu },{u_\mu },{\Theta _\mu },B}
 \right]} \right\|_N}\left( {\left\| {\nabla {\rho _\mu }} \right\|_{N - 1}^2 + \left\|
 {\left[ {{u_\mu },{\Theta _\mu }} \right]} \right\|_N^2} \right),\]
which will also be bounded by the right hand side term of
\eqref{2.16}.

Besides, for$|\alpha|\leq s$, standard energy estimates on
$\partial^\alpha E$ and $\partial^\alpha B$ from \eqref{2.2} yield
\begin{equation}\label{2.18}
\begin{split}
  &\frac{1}
{2}\frac{d}
{{dt}}\left( {{{\left\| {{\partial ^\alpha }E} \right\|}^2} + {{\left\| {{\partial ^\alpha }B}
 \right\|}^2}} \right) - \left\langle {\left( {1 + {\rho _e}} \right){\partial ^\alpha }{u_e} -
 \left( {1 + {\rho _i}} \right){\partial ^\alpha }{u_i},{\partial ^\alpha }E} \right\rangle  \\
   & = \left\langle {{\partial ^{\alpha  - \beta }}{\rho _e}{\partial ^\alpha }{u_e} - {\partial ^{\alpha
    - \beta }}{\rho _i}{\partial ^\alpha }{u_i},{\partial ^\alpha }E} \right\rangle \\
    & \leqslant
    C{\left\| E \right\|_s}\left( {\left\| {{u_\mu }} \right\|_s^2 + \left\| {\nabla {\rho _\mu }}
     \right\|_{s - 1}^2} \right),
\end{split}\end{equation}
which will be bounded by the right hand side term of \eqref{2.16}.
Then, with the help of \eqref{2.14}, the summation \eqref{2.17} and
\eqref{2.18} over $|\alpha|\leq s$, one has \eqref{2.16}.

\noindent \emph{Step 2.} It holds that
\begin{equation}\label{2.19}
\begin{split}
  \frac{d}
{{dt}} & \sum\limits_{\left| \alpha  \right| \leqslant s - 1}
{\sum\limits_{\mu  = e,i} {\left\langle {{\partial ^\alpha }{u_\mu
},\nabla {\partial ^\alpha }{\rho _\mu }} \right\rangle } }
 + \gamma \left( {\left\| {\nabla \left[ {{\rho _e},{\rho _i}} \right]} \right\|_{s- 1}^2
 + {{\left\| {{\rho _e} - {\rho _i}} \right\|}^2}} \right) \\
   & \leqslant C {\left( {\left\| {{u_\mu }} \right\|_s^2
   + \left\| {\left[ {{\rho _\mu },{u_\mu },{\Theta _\mu },B} \right]} \right\|_s^2\left( {\left\|
    {\nabla {\rho _\mu }} \right\|_{s - 1}^2 + \left\| {\left[ {{u_\mu },{\Theta _\mu }} \right]} \right\|_s^2}
     \right)} \right)} .
\end{split}\end{equation}
In fact, we can rewrite the equations \eqref{2.2} as
\begin{equation}
\label{2.20}
\left\{\begin{aligned}
&\partial_t {\rho_e} +\nabla\cdot u_e  =g_{1e}, \\
&\partial_t u_e +\nabla\rho_e+\nabla\Theta_e+ u_e+E=g_{2e}, \\
&\partial_t {\Theta_e} +\nabla\cdot u_e+ \Theta_e=g_{3e},
\\
&\partial_t {\rho_i} +\nabla\cdot u_i  =g_{1i},  \\
&\partial_t u_i +\nabla\rho_i+\nabla\Theta_i+u_i-E=g_{2i}
,\\
&\partial_t {\Theta_i} +\nabla\cdot u_i+ \Theta_i=g_{3i},\\
&  \partial_t E-\nabla\times B -u_e+ u_i= g_{4e}-g_{4i} ,
 \\
&  \partial_t B+\nabla\times E=0,\\
&     \nabla\cdot E=\rho_i-\rho_e ,\quad \nabla\cdot B=0,\quad
(t,x)\in(0,\infty)\times\mathbb{R}^3,
\end{aligned} \right.
\end{equation}
where
\begin{equation}
\label{2.21}
\left\{\begin{aligned}
&g_{1e}=-\rho_e \nabla\cdot u_e -u_e \nabla\rho_e, \\
&g_{2e}= -(u_e\cdot\nabla) u_e
-(\frac{\Theta_e+1}{1+\rho_e}-1)\nabla \rho_e- u_e \times B
, \\
&g_{3e}=-\Theta_e\nabla\cdot u_e-u_e\nabla\Theta_e,
\\
&g_{4e}=\rho_e
u_e,  \\
&g_{1i}=-\rho_i \nabla\cdot u_i -u_i \nabla\rho_i,  \\
&g_{2i}= -(u_i\cdot\nabla) u_i
-(\frac{\Theta_i+1}{1+\rho_i}-1)\nabla \rho_i+ u_i \times B
,\\
&g_{3i}=-\Theta_i\nabla\cdot u_i-u_i\nabla\Theta_i,\\
&g_{4i}= \rho_i u_i.
\end{aligned} \right.
\end{equation}
Let $\left| \alpha  \right| \leqslant s - 1$. Utilizing
$\partial^\alpha$ to the second equation of \eqref{2.20},
multiplying it by $\nabla\partial^\alpha\rho_e$, integrating over
$\mathbb{R}^3$ and using the last equation in \eqref{2.2}, replacing
$\partial_t {\rho_e}$ from the first equation of \eqref{2.20}
implies
\[\begin{gathered}
  \frac{d}
{{dt}}\left\langle {{\partial ^\alpha }{u_e},\nabla {\partial ^\alpha }{\rho _e}} \right\rangle
 + {\left\| {\nabla {\partial ^\alpha }{\rho _e}} \right\|^2} + {\left\| {{\partial ^\alpha }{\rho _e}}
  \right\|^2} - \left\langle {{\partial ^\alpha }{\rho _i},{\partial ^\alpha }{\rho _e}} \right\rangle
  + \left\langle {\nabla {\partial ^\alpha }{\Theta _e},\nabla {\partial ^\alpha }{\rho _e}} \right\rangle  \hfill \\
   = {\left\| {{\partial ^\alpha }\nabla  \cdot {u_e}} \right\|^2} + \left\langle {{\partial ^\alpha }
   \nabla {\rho _e},{\partial ^\alpha }{g_{2e}}} \right\rangle  - \left\langle {{\partial ^\alpha }{u_e},
   \nabla {\partial ^\alpha }{\rho _e}} \right\rangle  - \left\langle {{\partial ^\alpha }\nabla  \cdot {u_e},
   {\partial ^\alpha }{g_{1e}}} \right\rangle.  \hfill \\
\end{gathered} \]
Similarly as before, from the fourth and fifth equations of
\eqref{2.20}, we have
\[\begin{gathered}
  \frac{d}
{{dt}}\left\langle {{\partial ^\alpha }{u_i},\nabla {\partial ^\alpha }{\rho _i}} \right\rangle
+ {\left\| {\nabla {\partial ^\alpha }{\rho _i}} \right\|^2} + {\left\| {{\partial ^\alpha }{\rho _i}}
\right\|^2} - \left\langle {{\partial ^\alpha }{\rho _i},{\partial ^\alpha }{\rho _e}} \right\rangle
+ \left\langle {\nabla {\partial ^\alpha }{\Theta _i},\nabla {\partial ^\alpha }{\rho _i}} \right\rangle  \hfill \\
   = {\left\| {{\partial ^\alpha }\nabla  \cdot {u_i}} \right\|^2} + \left\langle {{\partial ^\alpha }\nabla
    {\rho _i},{\partial ^\alpha }{g_{2i}}} \right\rangle  - \left\langle {{\partial ^\alpha }{u_i},\nabla
     {\partial ^\alpha }{\rho _i}} \right\rangle  - \left\langle {{\partial ^\alpha }\nabla  \cdot {u_i},
     {\partial ^\alpha }{g_{1i}}} \right\rangle.  \hfill \\
\end{gathered} \]
Furthermore, the summation of the two equations above gives
\begin{equation}\notag
\begin{split}
  \frac{d}
{{dt}}&\left( {\left\langle {{\partial ^\alpha }{u_e}, \nabla
{\partial ^\alpha }{\rho _e}} \right\rangle  + \left\langle
{{\partial ^\alpha }{u_i},\nabla {\partial ^\alpha }{\rho _i}}
 \right\rangle } \right) + {\left\| {\nabla {\partial ^\alpha }{\rho _e}} \right\|^2}
 + {\left\| {\nabla {\partial ^\alpha }{\rho _i}} \right\|^2} + {\left\| {{\partial ^\alpha }
 \left( {{\rho _e} - {\rho _i}} \right)} \right\|^2} \\
   &= {\left\| {{\partial ^\alpha }\nabla  \cdot {u_e}} \right\|^2} + {\left\| {{\partial ^\alpha }
   \nabla  \cdot {u_i}} \right\|^2} - \left\langle {\nabla {\partial ^\alpha }{\Theta _i},
   \nabla {\partial ^\alpha }{\rho _i}} \right\rangle  - \left\langle {\nabla {\partial ^\alpha }
   {\Theta _e},\nabla {\partial ^\alpha }{\rho _e}} \right\rangle \\
   &\quad + \left\langle {{\partial ^\alpha }
   \nabla {\rho _e},{\partial ^\alpha }{g_{2e}}} \right\rangle  - \left\langle {{\partial ^\alpha }{u_e},
   \nabla {\partial ^\alpha }{\rho _e}} \right\rangle
    - \left\langle {{\partial ^\alpha }\nabla  \cdot {u_e},{\partial ^\alpha }{g_{1e}}}
    \right\rangle\\
    &\quad+ \left\langle {{\partial ^\alpha }\nabla {\rho _i},{\partial ^\alpha }{g_{2i}}} \right\rangle
   - \left\langle {{\partial ^\alpha }{u_i},\nabla {\partial ^\alpha }{\rho _i}} \right\rangle
   - \left\langle {{\partial ^\alpha }\nabla  \cdot {u_i},{\partial ^\alpha }{g_{1i}}}
   \right\rangle.
\end{split}\end{equation}
Therefore, after using Cauchy-Schwarz inequality, one has
\begin{equation}\label{2.22}
\begin{split}
  \frac{d}
{{dt}}&\left( {\left\langle {{\partial ^\alpha }{u_e},\nabla
{\partial ^\alpha } {\rho _e}} \right\rangle  + \left\langle
{{\partial ^\alpha }{u_i},\nabla {\partial ^\alpha } {\rho _i}}
\right\rangle } \right) + \lambda \left( {{{\left\| {\nabla
{\partial ^\alpha }{\rho _e}} \right\|}^2} + {{\left\| {\nabla
{\partial ^\alpha }{\rho _i}} \right\|}^2} + {{\left\|
{{\partial ^\alpha }\left( {{\rho _e} - {\rho _i}} \right)} \right\|}^2}} \right)  \\
   &\leqslant C {\left( {{{\left\| {{\partial ^\alpha }\nabla
   \cdot {u_\mu }} \right\|}^2} + {{\left\| {{\partial ^\alpha }{u_\mu }} \right\|}^2}
   + {{\left\| {{\partial ^\alpha }\nabla {\Theta _\mu }} \right\|}^2} + {{\left\|
    {{\partial ^\alpha }{g_{1\mu }}} \right\|}^2} + {{\left\| {{\partial ^\alpha }{g_{2\mu }}}
     \right\|}^2}} \right)}.
\end{split}\end{equation}
From the definition of $g_{j\mu}$, $(j=1,2)$, one can check that
\begin{equation}\notag
{\left\| {{\partial ^\alpha }{g_{1\mu }}} \right\|^2} + {\left\|
{{\partial ^\alpha }{g_{2\mu }}} \right\|^2} \leqslant C\left\|
{\left[ {{\rho _\mu },{u_\mu },{\Theta _\mu },B} \right]}
\right\|_s^2\left( {\left\| {\nabla {\rho _\mu }} \right\|_{s - 1}^2
+ \left\| {{u_\mu }} \right\|_s^2 + \left\| {{\Theta _\mu }}
\right\|_s^2} \right),
 \end{equation}
Putting this into \eqref{2.22}, then, \eqref{2.19} follows by taking
summation over $\left| \alpha \right| \leqslant s - 1$.

\noindent \emph{Step 3.} It holds that
\begin{equation}\label{2.23}
  \begin{split}
  \frac{d}
{{dt}}\sum\limits_{\left| \alpha  \right| \leqslant s - 1}
{\left\langle {{\partial ^\alpha }\left( {{u_e}
 - {u_i}} \right),{\partial ^\alpha }E} \right\rangle }  + \gamma \left\| E \right\|_{s - 1}^2
   \leqslant& C   \left\| {\left[ {{u_\mu },{\Theta _\mu }} \right]}
  \right\|_s^2 + C \|\nabla \rho_\mu \|_{s-1}^2+ C \left\| {{u_\mu }}
   \right\|_s \\
  \cdot{{\left\| {\nabla B} \right\|}_{s - 2}}
    &+ C \left\| U \right\|_s^2\left( {\left\|
   {\nabla {\rho _\mu }} \right\|_{s - 1}^2 + \left\| {\left[ {{u_\mu },{\Theta _\mu }} \right]} \right\|_s^2}
    \right) .
        \end{split}
\end{equation}
In fact, for $|\alpha|\leq s-1$, from the second and fifth equation
of \eqref{2.20}, one has
\begin{equation}\label{2.24}
  \begin{split}{\partial _t}\left( {{u_e} - {u_i}} \right) + \nabla \left( {{\rho
_e} - {\rho _i}} \right) + \nabla \left( {{\Theta _e} - {\Theta _i}}
\right) + 2E = {g_{2e}} - {g_{2i}} - \left( {{u_e} - {u_i}} \right).
\end{split}
\end{equation}
Utilizing $\partial^\alpha$ to \eqref{2.24},  multiplying it by
$\partial^\alpha E$, integrating over $\mathbb{R}^3$  and replacing
$\partial_t E$ from the seventh equation of \eqref{2.2} implies
\begin{equation}\notag
  \begin{split}
  \frac{d}
{{dt}}&\left\langle {{\partial ^\alpha }\left( {{u_e} - {u_i}}
\right),{\partial ^\alpha }E} \right\rangle + {\left\| {{\partial
^\alpha }\left( {{\rho _e} - {\rho _i}} \right)} \right\|^2} +
2{\left\| {{\partial ^\alpha }E}
\right\|^2}  \\
   &=  - \left\langle {{\partial ^\alpha }\left( {{\Theta _e} - {\Theta _i}} \right),{\partial ^\alpha }\left( {{\rho _e}
    - {\rho _i}} \right)} \right\rangle  + \left\langle {{\partial ^\alpha }\left( {{u_e} - {u_i}} \right),
    {\partial ^\alpha }E} \right\rangle  + \left\langle {{\partial ^\alpha }\left( {{u_e} - {u_i}} \right),
    \nabla  \times {\partial ^\alpha }B} \right\rangle   \\
   &\quad+ {\left\| {{\partial ^\alpha }\left( {{u_e} - {u_i}} \right)} \right\|^2} + \left\langle {{\partial ^\alpha }
   \left( {{u_e} - {u_i}} \right),{\partial ^\alpha }\left( {{\rho _e}{u_e} - {\rho _i}{u_i}} \right)} \right\rangle
   + \left\langle {{\partial ^\alpha }\left( {{g_{2e}} - {g_{2i}}} \right),{\partial ^\alpha }E}
   \right\rangle,
\end{split}
\end{equation}
Therefore, after using Cauchy-Schwarz inequality, one has
\begin{equation}\notag
  \begin{split}
  \frac{d}
{{dt}} & \left\langle {{\partial ^\alpha }\left( {{u_e} - {u_i}}
\right),{\partial ^\alpha }E}
 \right\rangle  + \gamma {\left\| {{\partial ^\alpha }E} \right\|^2}  \\
   \leq & C\left( {{{\left\| {{\partial ^\alpha } {{u_\mu} } } \right\|}^2}
   + {{\left\| {{\partial ^\alpha }{{\Theta _\mu}  }} \right\|}^2}}+
   {{\left\| {{\partial ^\alpha }{{\nabla\rho _\mu}  }} \right\|}^2}
   \right) + C{\left\| {\left[ {{u_e},{u_i}} \right]} \right\|_s}{\left\| {\nabla B} \right\|_{s - 2}}  \\
   &  + C\left\| {\left[ {{\rho _\mu },{u_\mu },{\Theta _\mu },B} \right]} \right\|_s^2\left( {\left\|
    {\nabla {\rho _\mu }} \right\|_{s - 1}^2 + \left\| {\left[ {{u_\mu },{\Theta _\mu }} \right]}
    \right\|_s^2} \right).
\end{split}
\end{equation}
Thus, with help of the summation of the previous estimate over
$|\alpha|\leq s-1$, one can obtain \eqref{2.23}.

\noindent \emph{Step 4.} It holds that
\begin{equation}\label{2.25}
\begin{split}
\frac{d} {{dt}}\sum\limits_{\left| \alpha  \right| \leqslant s - 2}
{\left\langle {{\partial ^\alpha }E, - \nabla  \times {\partial
^\alpha }B} \right\rangle  + \gamma \left\| {\nabla B} \right\|_{s -
2}^2}
 \leqslant  C   {(\left\| {\left[
{{u_\mu },E} \right]} \right\|_{s - 1}^2 + } \left\| {\nabla {\rho
_\mu }} \right\|_{s - 1}^2\left\| {{u_\mu }} \right\|_s^2 ).
\end{split}\end{equation}

In fact, for $|\alpha|\leq s - 2$, applying $\partial^\alpha$ to the
seventh equation of \eqref{2.2}, multiplying it by
$-\partial^\alpha\nabla\times B $, integrating over $\mathbb{R}^3$
and then utilizing the eighth equation of \eqref{2.2} gives
\begin{equation}\notag
\begin{split}
  \frac{d}
{{dt}}&\sum\limits_{\left| \alpha  \right| \leqslant s - 2}
{\left\langle {{\partial ^\alpha }E, - \nabla
 \times {\partial ^\alpha }B} \right\rangle  + {{\left\| {\nabla  \times {\partial ^\alpha }B} \right\|}^2}}   \\
   &= {\left\| {\nabla  \times {\partial ^\alpha }E} \right\|^2} - \left\langle {{\partial ^\alpha }\left( {{u_e}
    - {u_i}} \right),\nabla  \times {\partial ^\alpha }B} \right\rangle  + \left\langle {{\partial ^\alpha }
    \left( {{\rho _e}{u_e} - {\rho _i}{u_i}} \right), - \nabla  \times {\partial ^\alpha }B} \right\rangle
\end{split}
\end{equation}
Furthermore, with the help of Cauchy-Schwarz inequality and the
summation over $|\alpha|\leq s - 2$,  we yield \eqref{2.25}. Where
we have used
$$\left\|{\partial^\alpha\partial_i B}\right\|= \left\|{\partial_i\triangle^{-1}\nabla\times(\nabla\times
\partial^\alpha B)}\right\|
\leq C\left\|{\nabla\times\partial^\alpha B}\right\|$$ for $1\leq
i\leq 3$, due to $\nabla\cdot B=0$ and the fact that
$\partial_i\triangle^{-1}\nabla$ is bounded from $L^p$ to $L^p$ with
$1<p<\infty$, see \cite{Stein}.\\

\emph{Step 5.} Now, based on the four previous steps, we will search
\eqref{2.15}. We define the energy functional as
\begin{equation}\notag
 \begin{split}
 \mathcal{E}_s(U(t))=& \left\| U \right\|_s^2 +\mathcal {K}_1 \sum\limits_{\left| \alpha  \right| \leqslant s - 1}
  {\sum\limits_{\mu  = e,i}
   {\left\langle {{\partial ^\alpha }{u_\mu },\nabla {\partial ^\alpha }{\rho _\mu }} \right\rangle } }   \\
    &
+\mathcal {K}_2 \sum\limits_{\left| \alpha  \right| \leqslant s - 1}
{\left\langle {{\partial ^\alpha }\left( {{u_e}
    - {u_i}} \right),{\partial ^\alpha }E} \right\rangle }
  +\mathcal {K}_3 \sum\limits_{\left| \alpha  \right| \leqslant s - 2}
  {\left\langle {{\partial ^\alpha }E, - \nabla  \times {\partial ^\alpha }B} \right\rangle } ,
  \end{split}
\end{equation}
for constants $0<\mathcal {K}_3\ll\mathcal {K}_2\ll\mathcal {K}_1\ll
1$ to be chosen later. Notice that as soon as $0<\mathcal {K}_j
\ll1$ $_{(1\leq j \leq 3)}$ is sufficiently small, then $\mathcal
{E}_s(U(t))\sim ||U||^2_s$ holds true. Furthermore, the summation of
\eqref{2.16}, \eqref{2.19}$\times\mathcal {K}_1$,
\eqref{2.23}$\times\mathcal {K}_2$ and \eqref{2.25}$\times\mathcal
{K}_3$ implies that there is $0<\gamma<1$ such that
\begin{equation}\notag
\begin{split}
\frac{d}{dt}&\mathcal
{E}_s(U(t))+\|[u_e,u_i,\Theta_e,\Theta_i]\|_s^2 + \gamma\mathcal
{K}_1 (\|\nabla[\rho_e,\rho_i]\|_{s-1}^2+\|\rho_e-\rho_i\|^2)\\
&\quad+ \gamma\mathcal {K}_2\|E\|_{s-1}^2+ \gamma\mathcal {K}_3
\left\| {\nabla B}
\right\|_{s - 2}^2 \\
&\leq
    C[\mathcal {E}_s(U(t))^\frac{1}{2}+\mathcal {E}_s(U(t))]\mathcal
{D}_s(U(t))+C \mathcal {K}_1 \|u_\mu\|_s^2+C\mathcal
{K}_2\left(\|[u_\mu,\Theta_\mu]\|_s^2+\|\nabla\rho_\mu\|_{s-1}
^2\right)\\
&\quad  +C
\mathcal {K}_2 \|u_\mu\|_s\|\nabla B\|_{s-2} +C \mathcal {K}_3 \left\| {[u_\mu,E]} \right\|_{s - 1}^2\\
 &\leq
C[\mathcal {E}_s(U(t))^\frac{1}{2}+\mathcal {E}_s(U(t))]\mathcal
{D}_s(U(t))+C \mathcal {K}_1 \|u_\mu\|_s^2+C\mathcal
{K}_2\left(\|[u_\mu,\Theta_\mu]\|_s^2+\|\nabla\rho_\mu\|_{s-1}
^2\right)\\
&\quad   + \frac{1}{2}C\left(\mathcal {K}_2^\frac{1}{2}\left\|{u_\mu
}\right\|_s^2+\mathcal {K}_2^\frac{3}{2}\left\|{\nabla
B}\right\|_{s-2}^2\right) +C \mathcal {K}_3 \left\| {[u_\mu,E]}
\right\|_{s - 1}^2.
\end{split}
\end{equation}
 By letting
$0<\mathcal {K}_3\ll\mathcal {K}_2\ll\mathcal {K}_1\ll 1$ be
sufficiently small with $\mathcal {K}_2^{\frac{3}{2}}\ll\mathcal
{K}_3$, we obtain \eqref{2.15}. Now, we complete the proof of the
Theorem \ref{thm2.1}. \hfill $\Box$ \vspace{0.2cm}
%
\section{Linearized homogeneous equations}
In this section, for searching the time-decay property of solutions
to the nonlinear equations \eqref{2.2} in the last section, we have
to consider the decay properties of the linearized equations
\eqref{2.20}. Let us introduce the transformation
\begin{equation}\label{3.1}
\rho_1=\frac{\rho_e - \rho_i}{2},~u_1=\frac{u_e -
u_i}{2},~\Theta_1=\frac{\Theta_e - \Theta_i}{2}.
\end{equation}
Then, from system \eqref{2.2}, $U_1=[\rho_1,~u_1,~\Theta_1,~E,~B]$
satisfies
\begin{equation}\label{3.2}
\left\{\begin{split}
  &{\partial _t}{\rho _1} + \nabla  \cdot {u_1} = \frac{1}
{2}\left( {{g_{1e}} - {g_{1i}}} \right), \\
 & {\partial _t}{u_1} + \nabla {\rho _1} + \nabla {\Theta _1} + E + {u_1} = \frac{1}
{2}\left( {{g_{2e}} - {g_{2i}}} \right), \\
  &{\partial _t}{\Theta _1} + \nabla  \cdot {u_1} + {\Theta _1} = \frac{1}
{2}\left( {{g_{3e}} - {g_{3i}}} \right),  \\
 & {\partial _t}E - \nabla  \times B - 2{u_1} = {g_{4e}} - {g_{4i}}, \\
 & {\partial _t}B + \nabla  \times E = 0,  \\
 & \frac{1}
{2}\nabla  \cdot E =  - {\rho _1},\nabla  \cdot B = 0,\quad
(t,x)\in(0,\infty)\times\mathbb{R}^3,
\end{split}\right. \end{equation}
with initial value ${U_1}{|_{t = 0}} = {U_{1,0}}: = \left[ {{\rho
_{1,0}},{u_{1,0}},{\Theta _{1,0}},{E_0},{B_0}} \right],x \in
{\mathbb{R}^3}$ which satisfies the compatibility conditions
$\displaystyle\frac{1} {2}\nabla  \cdot E_0 =  - {\rho
_{1,0}},\nabla \cdot B_0 = 0.$ Where, $\left[ {{\rho
_{1,0}},{u_{1,0}},{\Theta _{1,0}}} \right]$ is given from $\left[
{{\rho _{\mu0}},{u_{\mu0}},{\Theta _{\mu0}}} \right]$ from the
transformation \eqref{3.1}.
Moreover, we introduce another transformation
\begin{equation}\label{3.3}
\rho_2=\frac{\rho_e + \rho_i}{2},~u_2=\frac{u_e +
u_i}{2},~\Theta_2=\frac{\Theta_e + \Theta_i}{2}.
\end{equation}
Then $U_2=[\rho_2,~u_2,~\Theta_2]$ satisfies
\begin{equation}\label{3.4}
\left\{\begin{split}
 & {\partial _t}{\rho _2} + \nabla  \cdot {u_2} = \frac{1}
{2}\left( {{g_{1e}} + {g_{1i}}} \right), \\
 & {\partial _t}{u_2} + \nabla {\rho _2} + \nabla {\Theta _2} + {u_2} = \frac{1}
{2}\left( {{g_{2e}} + {g_{2i}}} \right),\\
&  {\partial _t}{\Theta _2} + \nabla  \cdot {u_2} + {\Theta _2} =
\frac{1} {2}\left( {{g_{3e}} + {g_{3i}}} \right),\quad
(t,x)\in(0,\infty)\times\mathbb{R}^3,
\end{split}\right. \end{equation}
with  initial  value ${U_2}{|_{t = 0}} = {U_{2,0}}: = \left[ {{\rho
_{2,0}},{u_{2,0}},{\Theta _{2,0}}} \right],x \in {\mathbb{R}^3}$,
where $\left[ {{\rho _{2,0}},{u_{2,0}},{\Theta _{2,0}}} \right]$ is
from the transformation \eqref{3.3}.
Therefore, one can define the solution
$U_1=[\rho_1,~u_1,~\Theta_1,~E,~B]$ and
$U_2=[\rho_2,~u_2,~\Theta_2]$, respectively, as follows
\begin{equation}\label{3.5}
  {U_1}(t) = {e^{tL_1}}{U_{1,0}} + \frac{1}
{2}\int_0^t {{e^{\left( {t - y} \right)L_1}}\left[ {{g_{1e}} -
{g_{1i}},{g_{2e}} - {g_{2i}},{g_{3e}}
 - {g_{3i}},2\left( {{g_{4e}} - {g_{4i}}} \right)} \right](y)} dy,
\end{equation}
and
\begin{equation}\label{3.6}
  {U_2}(t) = {e^{tL_2}}{U_{2,0}} + \frac{1}
{2}\int_0^t {{e^{\left( {t - y} \right)L_2}}\left[ {{g_{1e}} +
{g_{1i}},{g_{2e}} + {g_{2i}},{g_{3e}}
 + {g_{3i}}} \right](y)} dy,
\end{equation}
where $e^{tL_1}U_{1,0}$  and $e^{tL_2}U_{2,0}$, respectively, denote
the solution of the following hohomogeneous initial value problems
\eqref{3.7}-\eqref{3.8} and \eqref{3.10}-\eqref{3.11}, which will be
given as follows:

The linearized homogeneous equations corresponding to \eqref{3.2} is
\begin{equation}\label{3.7}
\left\{\begin{split}
  &{\partial _t}{\rho _1} + \nabla  \cdot {u_1} = 0, \\
 & {\partial _t}{u_1} + \nabla {\rho _1} + \nabla {\Theta _1} + E + {u_1} = 0, \\
  &{\partial _t}{\Theta _1} + \nabla  \cdot {u_1} + {\Theta _1} = 0,  \\
 & {\partial _t}E - \nabla  \times B - 2{u_1} = 0, \\
 & {\partial _t}B + \nabla  \times E = 0,  \\
 & \frac{1}
{2}\nabla  \cdot E =  - {\rho _1},\nabla  \cdot B = 0,\quad
(t,x)\in(0,\infty)\times\mathbb{R}^3,
\end{split}\right. \end{equation}
with initial value
\begin{equation}\label{3.8}
{U_1}{|_{t = 0}} = {U_{1,0}}: = \left[ {{\rho
_{1,0}},{u_{1,0}},{\Theta _{1,0}},{E_0},{B_0}} \right],x \in
{\mathbb{R}^3}\end{equation}
  which satisfies the compatible
conditions
\begin{equation}\label{3.9}
\displaystyle\frac{1} {2}\nabla  \cdot E_0 =  - {\rho _{1,0}},\nabla
\cdot B_0 = 0.
\end{equation}
And the linearized homogeneous equations corresponding to
\eqref{3.7} is
\begin{equation}\label{3.10}
\left\{\begin{split}
 & {\partial _t}{\rho _2} + \nabla  \cdot {u_2} = 0, \\
 & {\partial _t}{u_2} + \nabla {\rho _2} + \nabla {\Theta _2} + {u_2} = 0,\\
&  {\partial _t}{\Theta _2} + \nabla  \cdot {u_2} + {\Theta _2}
=0,\quad (t,x)\in(0,\infty)\times\mathbb{R}^3,
\end{split}\right. \end{equation}
with initial value
\begin{equation}\label{3.11}
{U_2}{|_{t = 0}} = {U_{2,0}}: = \left[ {{\rho
_{2,0}},{u_{2,0}},{\Theta _{2,0}}} \right],x \in
{\mathbb{R}^3}.\end{equation}
Here $\left[ {{\rho _{2,0}},{u_{2,0}},{\Theta _{2,0}}} \right]$ is
from the transform \eqref{3.6}.
In the sequel, we usually denote $U_1$ $=[\rho_1$, $u_1$,
$\Theta_1$, $E$, $B]$ as the solution of the linearized homogeneous
equations \eqref{3.7}, and $U_2$ $=[\rho_2$, $u_2$, $\Theta_2]$ as
the one of \eqref{3.10}.

Firstly, for the linearized homogeneous system
\eqref{3.7}-\eqref{3.8}, similarly as \cite{FWK11}, we obtain the
$L^p-L^q$ decay property as follows

\begin{prop}\label{prop3.1}
Assume $U_1(t)=e^{tL_1}U_{1,0}$ is the solution to the initial value
problem \eqref{3.7}-\eqref{3.8} with ${U_{1,0}}$ $= [ {\rho
_{1,0}}$, ${u_{1,0}}$, ${\Theta _{1,0}}$, ${E_0}$, ${B_0} ]$ which
satisfies \eqref{3.8}. Then, for any $t\geq0$, $U_1$ $=[\rho_1$,
$u_1$, $\Theta_1$, $E$, $B]$ satisfies
\begin{equation}\label{3.12}
\left\{
\begin{split}
  &\left\| [{\rho_{1}(t),~ \Theta_{1} \left( t \right)}] \right\| \leqslant C{e^{ - \frac{t}
{2}}}\left\| {\left[ {{\rho _{1,0}},{u_{1,0}},{\Theta _{1,0}}} \right]} \right\|,   \\
 & \left\| {u_1\left( t \right)} \right\| \leqslant C{e^{ - \frac{t}
{2}}}\left\| {\left[ {{\rho _{1,0}},{\Theta _{1,0}}} \right]}
\right\| + C{\left( {1 + t} \right)^{ - \frac{5}
{4}}}{\left\| {\left[ {{u_{1,0}},{E_0},{B_0}} \right]} \right\|_{{L^1} \cap {{\dot H}^2}}},   \\
  & \left\| {E\left( t \right)} \right\| \leqslant C{\left( {1 + t} \right)^{ - \frac{5}
{4}}}{\left\| {\left[ {{u_{1,0}},{\Theta_{1,0}},{E_0},{B_0}} \right]} \right\|_{{L^1} \cap {{\dot H}^3}}},   \\
 & \left\| {B\left( t \right)} \right\| \leqslant C{\left( {1 + t} \right)^{ - \frac{3}
{4}}}{\left\| {\left[ {{u_{1,0}},{E_0},{B_0}} \right]}
\right\|_{{L^1} \cap {{\dot H}^2}}},
\end{split}\right.
\end{equation}
\begin{equation}\label{3.13}
\left\{\begin{split}& {\left\| [{\rho_{1}(t),~\Theta_{1} \left( t
\right)}] \right\|_{{L^\infty }}} \leqslant C{e^{ - \frac{t}
{2}}}{\left\| {\left[ {{\rho _{1,0}},{u_{1,0}},{\Theta _{1,0}}}
\right]} \right\|_{{L^2} \cap
{{\dot H}^2}}},\\
 &  {\left\| {u_1 \left( t \right)}
\right\|_{{L^\infty }}} \leqslant C{e^{ - \frac{t} {2}}}{\left\|
{\left[ {{\rho _{1,0}},{\Theta _{1,0}}} \right]} \right\|_{{L^1}
\cap {{\dot H}^2}}} + C{\left( {1 + t} \right)^{ - 2}}{\left\|
{\left[ {{u_{1,0}},{E_0},{B_0}} \right]} \right\|_{{L^1}
\cap {{\dot H}^5}}},  \\
&   {\left\| {E\left( t \right)} \right\|_{{L^\infty }}} \leqslant
C{\left( {1 + t} \right)^{ - 2}}
  {\left\| {\left[ {{u_{1,0}},{\Theta_{1,0}},{E_0},{B_0}} \right]} \right\|_{{L^1} \cap {{\dot H}^6}}}, \\
&   {\left\| {B\left( t \right)} \right\|_{{L^\infty }}} \leqslant
C{\left( {1 + t} \right)^{ - \frac{3} {2}}}{\left\| {\left[
{{u_{1,0}},{E_0},{B_0}} \right]} \right\|_{{L^1} \cap {{\dot
H}^5}}},
\end{split}\right.
\end{equation}
and
\begin{equation}\label{3.14}
\left\{\begin{split}
 & \left\| {\nabla B\left( t \right)} \right\| \leqslant C{\left( {1 + t} \right)^{ - \frac{5}
{4}}}{\left\| {\left[ {{u_{1,0}},{E_0},{B_0}} \right]} \right\|_{{L^1} \cap {{\dot H}^4}}},  \\
 & \left\| {{\nabla ^s}\left[ {E\left( t \right),B\left( t \right)} \right]} \right\| \leqslant
  C{\left( {1 + t} \right)^{ - \frac{5}
{4}}}{\left\| {\left[ {{u_{1,0}},\Theta_{1,0},{E_0},{B_0}} \right]}
\right\|_{{L^2} \cap {{\dot H}^{s + 3}}}}.
\end{split}\right.
\end{equation}
\end{prop}
%
\subsection{Explicit solutions of \eqref{3.10}-\eqref{3.11}}
Firstly, let us search the explicit Fourier transform solution
$U_2=[\rho_2,~u_2,~\Theta_2]$ of the initial value problem
\eqref{3.10}-\eqref{3.11}.

From the three equations of \eqref{3.10}, one has
\begin{equation}\label{3.15}
{\partial _{ttt}}{\rho _2} + 2{\partial _{tt}}{\rho _2} - 2\Delta {\partial _t}{\rho _2}
 + {\partial _t}{\rho _2} - \Delta {\rho _2} = 0,\end{equation}
with initial value
\begin{equation}\label{3.16}
\left\{
\begin{split}
&   {\rho _2}\left| {_{t = 0}} \right. = \rho _{2,0},  \\
  &{\partial _t}{\rho _2}\left| {_{t = 0}} \right. =  - \nabla  \cdot {u_{2,0}},  \\
   &{\partial _{tt}}{\rho _2}\left| {_{t = 0}} \right. = \Delta {\rho _{2,0}} + \nabla  \cdot {u_{2,0}} + \Delta {\Theta _{2,0}}.
  \end{split} \right.\end{equation}
After taking the Fourier transform on \eqref{3.15} and \eqref{3.16},
it follows that
\begin{equation}\label{3.17}
{\partial _{ttt}}{\hat{\rho} _2} + 2{\partial _{tt}}{\hat{\rho} _2}
 + (1+2|k|^2){\partial _t}{\hat{\rho} _2} +|k|^2 {\hat{\rho} _2} = 0,\end{equation}
with initial value
\begin{equation}\label{3.18}
\left\{
\begin{split}
&   {\hat{\rho} _2}\left| {_{t = 0}} \right. = \hat{\rho} _{2,0},  \\
  &{\partial _t}{\hat{\rho} _2}\left| {_{t = 0}} \right. =  - i|k|\tilde{k}  \cdot {\hat{u}_{2,0}},  \\
   &{\partial _{tt}}{\hat{\rho} _2}\left| {_{t = 0}} \right. = -|k|^2 {\hat{\rho} _{2,0}} + i|k|\tilde{k}  \cdot {\hat{u}_{2,0}}
   -|k|^2{\hat{\Theta} _{2,0}},
  \end{split} \right.\end{equation}
in this paper, we set $\tilde{k}=\frac{k}{|k|}.$ The characteristic
equation of \eqref{3.17} is
$$ F(\mathcal {X}):=\mathcal {X}^3+2\mathcal {X}^2+\left(1+2|k|^2\right)\mathcal {X}+|k|^2=0.$$
For the roots of the previous characteristic equation and their
properties, we obtain
\begin{lemma}\label{L3.1}
Assume $|k|\neq0.$ Then, $F(\mathcal {X})=0,$ $\mathcal
{X}\in\mathbb{C}$ has a real root
$\sigma=\sigma(|k|)\in(-\frac{1}{2},0)$ and two conjugate complex
roots $\mathcal {X}_\pm=\beta\pm i\omega$ with
$\beta=\beta(|k|)\in(-1,-\frac{3}{4})$ and
$\omega=\omega(|k|)\in(0,+\infty)$ which satisfy the following
properties:
\begin{equation}\label{3.19}
\beta=-1-\frac{\sigma}{2},~
\omega=\frac{1}{2}\sqrt{3\sigma^2+4\sigma+8|k|^2}.
\end{equation}
$\sigma,\beta,\omega$ are smooth in $|k|>0$, and $\sigma(|k|)$ is
strictly decreasing over $|k|>0$, with
$$ \lim_{|k|\longrightarrow0}\sigma(|k|)=0,~\lim_{|k|\longrightarrow\infty}\sigma(|k|)=-\frac{1}{2}.$$
Furthermore, the asymptotic behavior as follows hold true:
$$\sigma(|k|)=-O(1)|k|^2,~\beta(|k|)=-1+O(1)|k|^2,~\omega(|k|)=O(1)|k| $$
whenever $|k|\leq1$ is sufficiently small, and
$$\sigma(|k|)=-\frac{1}{2}+O(1)|k|^{-2},~\beta(|k|)=-\frac{3}{4}-O(1)|k|^{-2},~\omega(|k|)=O(1)|k|$$
whenever $|k|\geq1$ is sufficiently large. Here and in the sequel
$O(1)$ means strictly positive constant.
\end{lemma}

\noindent \emph{Proof.} Assume $|k|\neq0.$ First of all, we look for
the possibly existing real root for $F(\mathcal {X})=0$ over
$\mathcal {X}\in \mathbb{R}$. Since
 $$F'(\mathcal {X})=3\mathcal {X}^2+4\mathcal {X}+1+2|k|^{2}>0,$$
 and
 $F(-\frac{1}{2})=-\frac{1}{8}<0,~F(0)=|k|^{2}>0,$
 then equation $F(\mathcal {X})=0$ really has one and only one real
 root defined as $\sigma=\sigma(|k|)$ which satisfies $-\frac{1}{2}<\sigma<0.$
 After taking derivative of $F(\sigma(|k|))=0$ in $|k|$, one has
 $$\sigma'(|k|)=\frac{-|k|\left( 2+4\sigma\right)}{3\sigma^2+4\sigma+1+2|k|^{2}}<0,$$
 so that $\sigma(\cdot)$ is strictly decreasing in $|k|>0.$ Since
 $F(\sigma)=0$ can be re-written as
 $$\sigma\left[ \frac{\sigma(\sigma+2)}{1+2|k|^{2}}+1\right]=-\frac{|k|^{2}}{1+2|k|^{2}},$$
then $\sigma$ has limits $0$ and $-\frac{1}{2}$ as $|k|\rightarrow
0$ and $|k|\rightarrow \infty$, respectively.

$F(\sigma(|k|))=0$ is also equivalent with
$$\sigma+\frac{1}{2}=\frac{\frac{1}{2}(\sigma+1)^2}{(\sigma+1)^2+2|k|^{2}}$$
Therefore, it follows that $\sigma(|k|)=-O(1)|k|^2$ whenever $|k|<1$
is small enough and $\sigma(|k|)=-\frac{1}{2}+O(1)|k|^{-2} $
whenever $|k|\geq 1$ is large enough. Next, let us search roots of
$F(\mathcal {X})=0$ on $\mathcal {X}\in \mathbb{C}$. Since
$F(\sigma)=0$ with $\sigma\in \mathbb{R}$, $F(\mathcal {X})=0$ can
be split up into
$$ F(\mathcal {X})=(\mathcal {X}-\sigma)\left[\left(\mathcal {X}+1+\frac{\sigma}{2} \right)^2+\frac{3}{4}\sigma^2
+\sigma+2|k|^2 \right]=0.$$ Therefore, there are two conjugate
complex roots $\mathcal {X}_\pm=\beta\pm i\omega$ which satisfy
$$ \left(\mathcal {X}+1+\frac{\sigma}{2} \right)^2+\frac{3}{4}\sigma^2+\sigma+2|k|^2=0.$$
After solving the above equation, one can get that
$\beta=\beta(|k|),~\omega=\omega(|k|)$ take the form of
\eqref{3.19}. From the asymptotic behavior of $\sigma(|k|)$ at
$|k|=0$ and $\infty$, one can directly acquire that of
$\beta(|k|),~\omega(|k|)$. Now, we complete the proof of Lemma
\ref{L3.1}.\hfill$\Box$

Based on Lemma \ref{L3.1}, one can define the solution of
\eqref{3.17} as
\begin{equation}\label{3.20}
\hat{\rho}_2 (t,k)=c_1(k)e^{\sigma t}+e^{\beta t}\left(
c_2(k)\cos\omega t+c_3(k)\sin\omega t \right),
\end{equation}
where $c_i(k),~1\leq i \leq 3,$ is to be ascertained by \eqref{3.18}
later.
 In fact,
\eqref{3.18} implies
\begin{equation}\label{3.21}
\left[ {\begin{array}{*{20}{c}}
   {\hat \rho_2 {|_{t = 0}}}  \\
   {{\partial _t}\hat \rho_2 {|_{t = 0}}}  \\
   {{\partial _{tt}}\hat \rho_2 {|_{t = 0}}}  \\

 \end{array} } \right] = A\left[ {\begin{array}{*{20}{c}}
   {{c_1}}  \\
   {{c_2}}  \\
   {{c_3}}  \\

 \end{array} } \right],\quad A = \left[ {\begin{array}{*{20}{c}}
   1 & 1 & 0  \\
   \sigma  & \beta  & \omega   \\
   {{\sigma ^2}} & {{\beta ^2} - {\omega ^2}} & {2\beta \omega }  \\

 \end{array} } \right].\end{equation}
It is directly to check that
$$\det A = \omega \left[ {{\omega ^2} + {{\left( {\sigma  - \beta } \right)}^2}} \right]
= \omega \left( {3{\sigma ^2} + 4\sigma  + 1 + 2{{\left|
k \right|}^2}} \right) > 0
$$
and
\[{A^{ - 1}} = \frac{1}
{{\det A}}\left[ {\begin{array}{*{20}{c}}
   {\left( {{\beta ^2} + {\omega ^2}} \right)\omega } & { - 2\beta \omega } & \omega   \\
   {\sigma \left( {\sigma  - 2\beta } \right)\omega } & {2\beta \omega } & { - \omega }  \\
   {\sigma \left( {{\beta ^2} - {\omega ^2} - \sigma \beta } \right)} & {{\omega ^2} + {\sigma ^2}
   - {\beta ^2}} & {\beta  - \sigma }  \\

 \end{array} } \right].\]
Notice that \eqref{3.21} together with \eqref{3.18} gives
\begin{equation}\notag
\begin{split}
[{c_1},~ &   {c_2},~   {c_3}]^T =\frac{1} {{3{\sigma ^2} +
4\sigma  + 1 + 2 {{\left| k \right|}^2}}}\\
&\left[ {\begin{array}{*{20}{c}}
   {{\beta ^2} + {\omega ^2} -  { {{\left| k \right|}^2}} } & {  i\left| k \right|\left( {2\beta  + 1} \right)} & { - {{\left| k \right|}^2}}  \\
   {{\sigma ^2} - 2\sigma \beta  +  { {{\left| k \right|}^2}} } & {-i\left| k \right|\left( {2\beta  + 1} \right)} & {{{\left| k \right|}^2}}  \\
   {\frac{{\sigma \left( {{\beta ^2} - {\omega ^2} - \sigma \beta } \right) - \left( {\beta  - \sigma } \right) {{{\left| k \right|}^2}} }}
{\omega }} & {\frac{{i\left| k \right|}} {\omega }\left( {{\beta ^2}
- {\sigma ^2} - {\omega ^2} + \beta  - \sigma } \right)} &
{\frac{{\sigma  - \beta }}
{\omega }{{\left| k \right|}^2}}  \\
 \end{array} } \right]\left[ {\begin{array}{*{20}{c}}
   {{{\hat \rho }_{2,0}}}  \\
   {\tilde k \cdot {{\hat u}_{2,0}}}  \\
   {{{\hat \Theta }_{2,0}}}  \\
 \end{array} } \right].
 \end{split}\end{equation}
Here, we utilize $[\cdot]^T$  to denote the transpose of any vector.
Substituting the form of $\beta$ and $\omega$,  and making further
simplifications, we obtain
\begin{equation}\label{3.22}
\begin{split}
[&{c_1},~   {c_2},~   {c_3}]^T =\frac{1} {{3{\sigma ^2} +
4\sigma  + 1 + 2 {{\left| k \right|}^2}}}\\
&\left[ {\begin{array}{*{20}{c}}
   {{{\left( {\sigma  + 1} \right)}^2} + {{\left| k \right|}^2}} & { - i\left| k \right|\left( {\sigma  + 1} \right)} & { - {{\left| k \right|}^2}}  \\
   {2( \sigma  + 1) + {{\left| k \right|}^2}} & {i\left| k \right|\left( {\sigma  + 1} \right)} & {{{\left| k \right|}^2}}  \\
   {\frac{{{\sigma ( \sigma  + 1)} + (1- \frac{1}
{2}\sigma) {{\left| k \right|}^2}}} {\omega }} & {\frac{{i\left| k
\right|}} {\omega }\left( {\frac{3} {2}{\sigma ^2} + \frac{3}
{2}\sigma  + 2{{\left| k \right|}^2}} \right)} &
{\frac{{1 + \frac{3} {2}\sigma }}
{\omega }{{\left| k \right|}^2}}  \\
 \end{array} } \right]\left[ {\begin{array}{*{20}{c}}
   {{{\hat \rho }_{2,0}}}  \\
   {\tilde k \cdot {{\hat u}_{2,0}}}  \\
   {{{\hat \Theta }_{2,0}}}  \\
 \end{array} } \right].
 \end{split}\end{equation}
Similarly, from the three equations of \eqref{3.10}, one has
\begin{equation}\label{3.23}
\partial_{ttt}\hat{\Theta}_2
+2\partial_{tt}\hat{\Theta}_2+\left(1+2|k|^2\right)\partial_{t}\hat{\Theta}_2
+|k|^2\hat{\Theta}_2=0,
\end{equation}
with initial value
\begin{equation}\label{3.24}
\left\{\begin{split}
&\hat{\Theta}_2|_{t=0}=\hat{\Theta}_{2,0},\\
&\partial_{t}\hat{\Theta}_2|_{t=0}=-i|k|\tilde{k}\cdot \hat{u}_{2,0}-\hat{\Theta}_{2,0},\\
&\partial_{tt}\hat{\Theta}_2|_{t=0}=-|k|^2\hat{\rho}_{2,0}+2i|k|\tilde{k}\cdot
\hat{u}_{2,0}+\left(1-|k|^2\right)\hat{\Theta}_{2,0}.
\end{split}
\right.\end{equation} Based on Lemma \ref{L3.1}, one can also set
the solution of \eqref{3.23} as
\begin{equation}\label{3.25}
\hat{\Theta}_2(t,k)=c_4(k)e^{\sigma t}+e^{\beta t}\left(
c_5(k)\cos\omega t+c_6(k)\sin\omega t \right),
\end{equation}
where $c_i(k),~4\leq i \leq 6,$ is to be ascertained by \eqref{3.24}
later. In fact, after tenuous computation, \eqref{3.24} implies
\begin{equation}\label{3.26}
\begin{split}
[&{c_4},~   {c_5},~   {c_6}]^T =\frac{1} {{3{\sigma ^2} +
4\sigma  + 1 + 2 {{\left| k \right|}^2}}} \\
& \left[ {\begin{array}{*{20}{c}}
   {- {{\left| k \right|}^2}} & {-i\left| k
\right|\left( {1 + \sigma  } \right)} & {(1 + \sigma) \sigma+{{\left| k \right|}^2}}  \\
   {{{\left| k \right|}^2}} & { i\left| k
\right|\left( {1 + \sigma  } \right)} & {(1 + 2\sigma)(1 + \sigma) +{{\left| k \right|}^2} }  \\
   {\frac{\frac{3}{2}\sigma +1 } {\omega }|k|^2} &
{\frac{-i|k|} {\omega }(\frac{3}{2}\sigma(\sigma+2)+1+2|k|^2)} & {-\frac{|k|^2+\frac{1}{2}\sigma(|k|^2+1+\sigma)  }
{\omega }}  \\
 \end{array} } \right]\left[ {\begin{array}{*{20}{c}}
   {{{\hat \rho }_{2,0}}}  \\
   {\tilde k \cdot {{\hat u}_{2,0}}}  \\
   {{{\hat \Theta }_{2,0}}}  \\
 \end{array} } \right].
 \end{split}\end{equation}
Similarly, again from the three equations of \eqref{3.10}, we also
have
\begin{equation}\label{3.27}
{\partial _{ttt}} ({\tilde{k}\cdot \hat{u}_2}) + 2{\partial _{tt}}(\tilde{k}\cdot \hat{u}_2)+
(1+2|k|^2)\partial_t ({\tilde{k}\cdot \hat{u}_2})+|k|^2({\tilde{k}\cdot \hat{u}_2})= 0,
\end{equation}
with initial value
\begin{equation}\label{3.28}
\left\{\begin{split}
&{\tilde{k}\cdot \hat{u}_2}|_{t=0}={\tilde{k}\cdot \hat{u}}_{2,0},\\
&\partial_{t}{(\tilde{k}\cdot \hat{u}_2)}|_{t=0}=-i|k|{\hat
\rho _{2,0}} - \tilde k \cdot {\hat u_{2,0}} - i\left| k \right|{\hat \Theta
_{2,0}},
\\
&\partial_{tt}{(\tilde{k}\cdot \hat{u}_2)}|_{t=0}=i|k|{\hat \rho _{2,0}} +(1-2{\left| k \right|^2})\tilde k \cdot {\hat u_{2,0}} + 2i\left|
k \right|{\hat \Theta _{2,0}}.
\end{split}
\right.\end{equation}
 From Lemma \ref{L3.1}, one can also check that the
solution of \eqref{3.27} has the form
\begin{equation}\label{3.29}
\tilde{k}\cdot \hat{u}_2(t,k)=c_7(k)e^{\sigma t}+e^{\beta t}\left(
c_8(k)\cos\omega t+c_9(k)\sin\omega t \right),
\end{equation}
with
\begin{equation}\label{3.30}
 \begin{split} [&{c_7},~   {c_8},~   {c_9}]^T =\frac{1}
{{3{\sigma ^2} +
4\sigma  + 1 + 2 {{\left| k \right|}^2}}} \\
& \left[ {\begin{array}{*{20}{c}}
   {-i\left| k
\right|\left( {1 + \sigma  } \right)} & {\sigma\left( {1 + \sigma  } \right)  } & {-i\left| k
\right| \sigma  }  \\
   {i\left| k
\right|\left( {1 + \sigma  } \right)} & { \left( {1 + \sigma  } \right)\left( {1 + 2\sigma  } \right)+2|k|^2} & {i\left| k
\right| \sigma   }  \\
   {\frac{-i|k|} {\omega }(\frac{3}{2}\sigma(\sigma+1)-2|k|^2)} &
{   \frac{-\sigma \left( {1 + \sigma  } -2|k|^2 \right) }{2\omega}      } & {\frac{i|k|} {\omega }(-\frac{3}{2}\sigma (\sigma+2)+2|k|^2-1)  }
 \\
 \end{array} } \right]\left[ {\begin{array}{*{20}{c}}
   {{{\hat \rho }_{2,0}}}  \\
   {\tilde k \cdot {{\hat u}_{2,0}}}  \\
   {{{\hat \Theta }_{2,0}}}  \\
 \end{array} } \right].
 \end{split}\end{equation}

Furthermore, after taking the curl for the second equation of
\eqref{3.10} and making the Fourier transform in $x$, we have
\begin{equation}\label{3.31}
\partial_t \left(\tilde{k}\times  (\tilde{k}\times \hat{u}_2)\right)
+ \tilde{k}\times  (\tilde{k}\times \hat{u}_2)=0,
\end{equation}
with initial value
\begin{equation}\label{3.32}
\tilde{k}\times  (\tilde{k}\times \hat{u}_2)\left|_{t=0}\right.=
\tilde{k}\times (\tilde{k}\times \hat{u}_{2,0}).
\end{equation}
After solving \eqref{3.31}-\eqref{3.32}, we have
\begin{equation}\label{3.33}
\tilde{k}\times  (\tilde{k}\times \hat{u}_2)
=e^{-t}\left(\tilde{k}\times (\tilde{k}\times \hat{u}_{2,0})\right).
\end{equation}

 Now, we can obtain the explicit Fourier transform solution $U_2=[\rho_2,~u_2,~\Theta_2]$ as follows
 from the above computations.
\begin{theorem}\label{thm3.1}
Assume $U_2=[\rho_2,~u_2,~\Theta_2]$ be the solution of the initial
value problem \eqref{3.10}-\eqref{3.11} on the linearized
homogeneous equations. For $(t,k)\in(0,\infty)\times\mathbb{R}^3$
with $|k|\neq0$, we obtain
\begin{equation}\label{3.34}
\left[ {\begin{array}{*{20}c}
   {\hat{\rho }_2  (t,k)}  \\
   {\hat{u}_2  (t,k)}  \\
   {\hat{\Theta}_2 (t,k)}  \\
\end{array}} \right] = \left[ {\begin{array}{*{20}c}
   {\hat{\rho}_2  (t,k)}  \\
   {\hat{u}_{ 2||} (t,k)}  \\
    {\hat{\Theta}_2  (t,k)}  \\
   \end{array}} \right] + \left[ {\begin{array}{*{20}c}
   0  \\
   {\hat{u}_{ 2 \bot } (t,k)}  \\
          0\\
\end{array}} \right].
\end{equation}
Here $\hat{u}_{ 2||},\hat{u}_{  2\bot }$ are defined by
$$\hat{u}_{2
||}=\tilde{k}\tilde{k}\cdot\hat{u}_2,\ ~ \hat{u}_{ 2\bot
}=-\tilde{k}\times(\tilde{k}\times\hat{u}_2)=(I_3-\tilde{k}\otimes\tilde{k})
\hat{u}_2.$$
 Then, there exit matrices $G^I_{5\times5}(t,k)$ and
$G^{II}_{3\times3}(t,k)$ such that
\begin{equation}\label{3.35}
\left[ {\begin{array}{*{20}c}
   {\hat{\rho_2} (t,k)}  \\
      {\hat{u}_{2||} (t,k)}  \\
{\hat{\Theta}_2 (t,k)}  \\
   \end{array}} \right] = G_{5 \times 5}^I (t,k)\left[ {\begin{array}{*{20}c}
   {\hat{\rho} _{2,0} (k)}  \\
      {\hat{u}_{2||,0} (k)}  \\
{\hat{\Theta} _{2,0} (k)}  \\
   \end{array}} \right]
\end{equation}
and
\begin{equation}\label{3.36}
 {\begin{array}{*{20}c}
   {\hat{u}_{  2\bot } (t,k)}  \\
      \end{array}}  = G_{3 \times 3}^{II} (t,k) {\begin{array}{*{20}c}
   {\hat{u}_{  2\bot,0 } (k)}  \\
     \end{array}} ,
\end{equation}
where $ G_{5 \times 5}^I$ is explicitly ascertained by
representations \eqref{3.20}, \eqref{3.29}, \eqref{3.25} for
$\hat{\rho}_2(t,k)$, $\hat{u}_{2||}(t,k)$, $\hat{\Theta}_2(t,k)$
with $c_i(k)$, $(1\leq i\leq 9)$ are defined as \eqref{3.22},
\eqref{3.30}, \eqref{3.26} in terms of $\hat{\rho}_{2,0}(k)$,
$\hat{u}_{2||,0}(k)$, $\hat{\Theta}_{2,0}(k)$; and $G_{3 \times
3}^{II} $ is chosen by the representations \eqref{3.33} for
$\hat{u}_{ 2\bot } (t,k)$ in terms of $\hat{u}_{  2\bot,0 } (k)$.
\end{theorem}
\subsection{$L^p-L^q$ decay
property.} In this subsection, we use Theorem \ref{thm3.1} to obtain
$L^p-L^q$ decay property for every component of the solution $U_2$
$=[\rho_2$, $u_2$, $\Theta_2]$. For this aim, we first search the
rigorous time-frequency estimates on $\hat{U}_2$ $=[\hat{\rho}_2$,
$\hat{u}_2$, $\hat{\Theta}_2]$ as follws
\begin{lemma}\label{L3.2}
Assume $U_2=[\rho_2,~u_2,~\Theta_2]$ be the solution to the initial
value problem \eqref{3.10}-\eqref{3.11} on the linearized
homogeneous equations. Then, there are constants $\gamma>0,C>0$ such
that for all $(t,k)\in(0,\infty)\times\mathbb{R}^3$,
\begin{equation}\label{3.37}
\begin{split}
  \left|\hat{\rho}_2  (t,k) \right| \leqslant
   C\left| {\left[ {{{\hat \rho}_{2,0}}(t,k),{{\hat u}_{2,0}}(t,k),{{\hat \Theta}_{2,0}}(t,k)} \right]} \right| \cdot
    \left\{ {\begin{array}{*{20}{c}}
   {{e^{ - \gamma t}} + {e^{ - \gamma {{\left| k \right|}^2}t}}}
    & {if~\left| k \right| \leqslant 1,}  \\
   {{e^{ - \gamma t}} +
{e^{\frac{{ - \gamma }}
{{{{\left| k \right|}^2}}}t}}} & {if~\left| k \right| > 1,}  \\
 \end{array} } \right.
\end{split}
\end{equation}
\begin{equation}\label{3.38}
\begin{split}
  \left| {\hat u_2(t,k)} \right| \leqslant 
  C\left| {\left[ {{{\hat \rho}_{2,0}}(k),{{\hat u}_{2,0}}(k),{{\hat \Theta}_{2,0}}(k)} \right]} \right|
  \cdot
    \left\{ {\begin{array}{*{20}{c}}
   {{e^{ - \gamma t}} + \left| k \right|{e^{ - \gamma {{\left| k \right|}^2}t}}}
    & {if~\left| k \right| \leqslant 1,}  \\
   {{|k|^{-1}e^{ - \gamma t}} +
{e^{\frac{{ - \gamma }}
{{{{\left| k \right|}^2}}}t}}} & {if~\left| k \right| > 1,}  \\
 \end{array} } \right.
\end{split}
\end{equation}
and
\begin{equation}\label{3.39}
\begin{split}
  \left|\hat{\Theta}_2  (t,k) \right| \leqslant
   C\left| {\left[ {{{\hat \rho}_{2,0}}(t,k),{{\hat u}_{2,0}}(t,k),{{\hat \Theta}_{2,0}}(t,k)} \right]} \right| \cdot
    \left\{ {\begin{array}{*{20}{c}}
   {{e^{ - \gamma t}} + {e^{ - \gamma {{\left| k \right|}^2}t}}}
    & {if~\left| k \right| \leqslant 1,}  \\
   {{e^{ - \gamma t}} +
{e^{\frac{{ - \gamma }}
{{{{\left| k \right|}^2}}}t}}} & {if~\left| k \right| > 1,}  \\
 \end{array} } \right.
\end{split}
\end{equation}
\end{lemma}
\noindent \emph{Proof.}~  Firstly, let us look for the upper bound
of $\hat{\rho}_2$ defined as \eqref{3.37}.
 In fact, from Lemma \ref{L3.1}, it is directly to check
 \eqref{3.22} to get
\[\left[ {\begin{array}{*{20}{c}}
   {{c_1}}  \\
   {{c_2}}  \\
   {{c_3}}  \\
 \end{array} } \right] = \left[ {\begin{array}{*{20}{c}}
   {O(1){}} & { - O(1){{\left| k \right|}}i} & { - O(1){{\left| k \right|}^2}}  \\
   {O(1)} & {O(1){{\left| k \right|}}i} & {O(1){{\left| k \right|}^2}}  \\
   {O(1)\left| k \right|} & { - O(1)\left| k \right|^2i} & {  O(1){{\left| k \right|}}}  \\
 \end{array} } \right]\left[ {\begin{array}{*{20}{c}}
   {{{\hat \rho }_{2,0}}}  \\
   {\tilde k \cdot {{\hat u}_{2,0}}}  \\
   {{{\hat \Theta }_{2,0}}}  \\
 \end{array} } \right]\]
as $|k|\rightarrow0$, and
\[\left[ {\begin{array}{*{20}{c}}
   {{c_1}}  \\
   {{c_2}}  \\
   {{c_3}}  \\
 \end{array} } \right] = \left[ {\begin{array}{*{20}{c}}
   {O(1)} & { - O(1){{\left| k \right|}^{ - 1}}i} & { - O(1)}  \\
   {O(1)} & {O(1){{\left| k \right|}^{ - 1}}i} & {O(1)}  \\
   {O(1){{\left| k \right|}^{ - 1}}} & { - O(1)i} & {O(1){{\left| k \right|}^{ - 1}}}  \\
 \end{array} } \right]\left[ {\begin{array}{*{20}{c}}
   {{{\hat \rho }_{2,0}}}  \\
   {\tilde k \cdot {{\hat u}_{2,0}}}  \\
   {{{\hat \Theta }_{2,0}}}  \\
 \end{array} } \right]\]
as $|k|\rightarrow\infty$.

Therefore, after putting the previous computations into
\eqref{3.20}, it holds that
\[\begin{split}
  \hat \rho_2 (t,k) =& \left( {O(1){{}}{{\hat \rho }_{2,0}} - O(1)
  {{\left| k \right|}}i\tilde k \cdot {{\hat u}_{2,0}} - O(1){{\left| k \right|}^2}
  {{\hat \Theta }_{2,0}}} \right){e^{\sigma t}}   \\
   &+ \left( {O(1){{\hat \rho }_{2,0}} + O(1){{\left| k \right|}}i\tilde k \cdot {{\hat u}_{2,0}}
   + O(1){{\left| k \right|}^2}{{\hat \Theta }_{2,0}}} \right){e^{\beta t}}\cos \omega t  \\
   &+ \left( {O(1)|k|{{\hat \rho }_{2,0}} - O(1)\left| k \right|^2i\tilde k \cdot {{\hat u}_{2,0}} + O(1)
   {{\left| k \right|}}{{\hat \Theta }_{2,0}}} \right){e^{\beta t}}\sin \omega t,  \\
\end{split} \]
as $|k|\rightarrow0$, and
\[\begin{split}
  \hat \rho_2 (t,k) =& \left( {O(1){{\hat \rho }_{2,0}} - O(1){{\left| k \right|}^{ - 1}}i\tilde k \cdot
  {{\hat u}_{2,0}} - O(1){{\hat \Theta }_{2,0}}} \right){e^{\sigma t}}   \\
   &+ \left( {O(1){{\hat \rho }_{2,0}} + O(1){{\left| k \right|}^{ - 1}}i\tilde k \cdot {{\hat u}_{2,0}}
   + O(1){{\hat \Theta }_{2,0}}} \right){e^{\beta  t}}\cos \omega t  \\
   &+ \left( {O(1){{\left| k \right|}^{ - 1}}{{\hat \rho }_{2,0}} - O(1)i\tilde k \cdot {{\hat u}_{2,0}}
   + O(1){{\left| k \right|}^{ - 1}}{{\hat \Theta }_{2,0}}} \right){e^{\beta  t}}\sin \omega t, \\
\end{split} \]
as $|k|\rightarrow\infty$.

Based on Lemma \ref{L3.1}, we find that there is $\gamma>0$ such
that
\begin{equation}\notag\left\{\begin{split}
   & {\sigma (k) \leqslant  - \gamma {{\left| k \right|}^2},}\quad  {\beta (k) =  - 1 - \frac{\sigma }
{2} \leqslant  - \gamma }\quad  {{over}\left| k \right| \leqslant 1,}  \\
   &{\sigma (k) \leqslant  - \gamma ,} \quad {\beta (k) =  - 1 - \frac{\sigma }
{2} \leqslant  - \frac{\gamma } {{{{\left| k \right|}^2}}}}
\quad{{over}\left| k \right| \geqslant 1.}
  \end{split} \right.\end{equation}
Thus, one can obtain, for $|k|\leq 1$
\[\left| {{{\hat \rho }_2}\left( {t,k} \right)} \right| \leqslant
C\left( {{e^{ - \gamma t}} + {e^{ - \gamma {{\left| k
\right|}^2}t}}} \right)\left| {\left[ {{{\hat \rho }_{2,0}},\tilde k
\cdot {{\hat u}_{2,0}},{{\hat \Theta }_{2,0}}} \right]} \right|,\]
and for $|k|\geq 1$
\[\left| {{{\hat \rho }_2}\left( {t,k} \right)} \right| \leqslant
C\left( {{e^{ - \gamma t}} + {e^{ - \frac{\gamma}{{\left| k
\right|}^2} {}t}}} \right)\left| {\left[ {{{\hat \rho
}_{2,0}},\tilde k \cdot {{\hat u}_{2,0}},{{\hat \Theta }_{2,0}}}
\right]} \right|.\]
Furthermore, one has
\begin{equation}\notag
\left| {{{\hat \rho }_2}\left( {t,k} \right)} \right| \leqslant C
\left| {\left[ {{{\hat \rho }_{2,0}},  {{\hat u}_{2,0}},{{\hat
\Theta }_{2,0}}} \right]} \right|\cdot \left\{\begin{split}\left(
{{e^{ - \gamma t}} + {e^{ - \gamma {{\left| k
\right|}^2}t}}} \right) & \quad {if}\quad |k|\leq 1, \\
\left( {{e^{ - \gamma t}} + {e^{ - \frac{\gamma}{{\left| k
\right|}^2} {}t}}} \right) & \quad {if}\quad |k|\geq 1.
 \end{split} \right.\end{equation}
Similarly, we obtain \eqref{3.38} and \eqref{3.39}. Now, we complete
the proof of Lemma \ref{L3.2}. \hfill$\Box$

From Lemma \ref{L3.2}, it is straightforward to acquire the decay
property for every component of the solution $U_2$ $=[\rho_2$,
$u_2$, $\Theta_2]$. So that we omitted the details of proof for
briefness. See for instance \cite{FWK11}.
\begin{theorem}\label{thm3.2}
Assume $1\leq p,r\leq 2\leq q\leq\infty,l\geq0$ and an integer
$m\geq0$. Suppose
 $U_2(t)$ $=e^{tL_2}U_{2,0}$ be the solution of the initial value
problem \eqref{3.10}-\eqref{3.11}. Then, for any $t\geq0$, $U_2$
$=[\rho_2$, $u_2$, $\Theta_2]$ satisfies decay property as follows
\begin{equation}\label{3.40}
\begin{split}
  {\left\| {{\nabla ^m}{\rho _2}\left( t \right)} \right\|_{L^q}} \leqslant & C{\left( {1 + t} \right)^{ - \frac{3}
{2}\left( {\frac{1} {p} - \frac{1} {q}} \right) - \frac{m}
{2}}}{\left\| {\left[ {{\rho _{2,0}},{u_{2,0}},{\Theta _{2,0}}} \right]} \right\|_{{L^p}}}  \\
   &+ C{\left( {1 + t} \right)^{ - \frac{l}
{2}}}{\left\| {{\nabla ^{m + {{\left[ {l + 3\left( {\frac{1} {r} -
\frac{1} {q}} \right)} \right]}_ + }}}\left[ {{\rho
_{2,0}},{u_{2,0}},{\Theta _{2,0}}} \right]} \right\|_{{L^r}}},
\end{split}\end{equation}
\begin{equation}\label{3.41}
\begin{split}
  {\left\| {{\nabla ^m}u_2\left( t \right)} \right\|_{{L^q}}} \leqslant& C{\left( {1 + t} \right)^{ - \frac{3}
{2}\left( {\frac{1} {p} - \frac{1} {q}} \right) - \frac{{m + 1}}
{2}}}{\left\| {\left[ {{\rho _{2,0}},{u_{2,0}},{\Theta _{2,0}}} \right]} \right\|_{{L^p}}}  \\
  & + C{\left( {1 + t} \right)^{ - \frac{{l }}
{2}}}{\left\| {{\nabla ^{m + {{\left[ {l + 3\left( {\frac{1} {r} -
\frac{1}
{q}} \right)} \right]}_ + }}}\left[ {\rho _{2,0},{u_{2,0}},{\Theta _{2,0}}} \right]} \right\|_{{L^r}}},  \\
\end{split}
\end{equation}
and
\begin{equation}\label{3.42}
\begin{split}
  {\left\| {{\nabla ^m}{\Theta _2}\left( t \right)} \right\|_{L^q}} \leqslant & C{\left( {1 + t} \right)^{ - \frac{3}
{2}\left( {\frac{1} {p} - \frac{1} {q}} \right) - \frac{m}
{2}}}{\left\| {\left[ {{\rho _{2,0}},{u_{2,0}},{\Theta _{2,0}}} \right]} \right\|_{{L^p}}}  \\
   &+ C{\left( {1 + t} \right)^{ - \frac{l}
{2}}}{\left\| {{\nabla ^{m + {{\left[ {l + 3\left( {\frac{1} {r} -
\frac{1} {q}} \right)} \right]}_ + }}}\left[ {{\rho
_{2,0}},{u_{2,0}},{\Theta _{2,0}}} \right]} \right\|_{{L^r}}},
\end{split}
\end{equation}
 where
\begin{equation}\notag
{[l + 3(\frac{1}{r} - \frac{1}{q})]_ + } = \left\{
{\begin{array}{*{20}{c}}
   l & { \mbox{if} ~~  r = q = 2\ \mbox{ and }{\rm{ }}l{\mbox{ is an integer,}}}  \\
   {{{[l + 3(\frac{1}{r} - \frac{1}{q})]}_ - } + 1} & \mbox{otherwise,}  \\
\end{array}} \right.
\end{equation}
where, we use $[\cdot]_-$ to denote the integer part of the
argument.
\end{theorem}
From Theorem \ref{thm3.2}, let us list some particular cases as
follows for later use.
\begin{coro}\label{Corollary3.1}
Let $U_2(t)=e^{tL_2}U_{2,0}$ be the solution of the initial value
problem \eqref{3.10}-\eqref{3.11}. Then, for any $t\geq0$, $U_2$
$=[\rho_2$, $u_2$, $\Theta_2]$ satisfies
\begin{equation}\label{3.43}
\left\{
\begin{split}
  &\left\| {{\rho _2}\left( t \right)} \right\| \leqslant C{\left( {1 + t} \right)^{ - \frac{3}
{4}}}{\left\| {\left[ {{\rho _{2,0}},{u_{2,0}},{\Theta _{2,0}}}
\right]} \right\|_{{L^1} \cap {{\dot H}^2}}},   \\
 & \left\| {{u _2}\left( t \right)} \right\| \leqslant C{\left( {1 + t} \right)^{ - \frac{5}
{4}}}{\left\| {\left[ {{\rho _{2,0}},{u_{2,0}},{\Theta _{2,0}}}
\right]} \right\|_{{L^1} \cap {{\dot H}^3}}}, \\
 & \left\| {{\Theta _2}\left( t \right)} \right\| \leqslant C{\left( {1 + t} \right)^{ - \frac{3}
{4}}}{\left\| {\left[ {{\rho _{2,0}},{u_{2,0}},{\Theta _{2,0}}}
\right]} \right\|_{{L^1} \cap {{\dot H}^2}}},
\end{split}\right.
\end{equation}
\begin{equation}\label{3.44}
\left\{\begin{split}& \left\| {\nabla {\rho _2}\left( t \right)}
\right\| \leqslant C{\left( {1 + t} \right)^{ - \frac{5}
{4}}}{\left\| {\left[ {{\rho _{2,0}},{u_{2,0}},{\Theta _{2,0}}}
\right]} \right\|_{{L^1} \cap {{\dot H}^4}}},
\\
 &  \left\| {\nabla {u _2}\left( t \right)}
\right\| \leqslant C{\left( {1 + t} \right)^{ - \frac{7}
{4}}}{\left\| {\left[ {{\rho _{2,0}},{u_{2,0}},{\Theta _{2,0}}}
\right]} \right\|_{{L^1} \cap {{\dot H}^5}}},  \\
&   \left\| {\nabla {\Theta _2}\left( t \right)} \right\| \leqslant
C{\left( {1 + t} \right)^{ - \frac{5} {4}}}{\left\| {\left[ {{\rho
_{2,0}},{u_{2,0}},{\Theta _{2,0}}} \right]} \right\|_{{L^1} \cap
{{\dot H}^4}}}
\end{split}\right.
\end{equation}
and
\begin{equation}\label{3.45}
\left\{
\begin{split}
  &\left\| {{\rho _2}\left( t \right)} \right\|_{L^\infty} \leqslant C{\left( {1 + t} \right)^{ - \frac{3}
{2}}}{\left\| {\left[ {{\rho _{2,0}},{u_{2,0}},{\Theta _{2,0}}}
\right]} \right\|_{{L^1} \cap {{\dot H}^5}}},   \\
 & \left\| {{u _2}\left( t \right)} \right\|_{L^\infty} \leqslant C{\left( {1 + t} \right)^{ - 2}}
 {\left\| {\left[ {{\rho _{2,0}},{u_{2,0}},{\Theta _{2,0}}}
\right]} \right\|_{{L^1} \cap {{\dot H}^6}}}, \\
 & \left\| {{\Theta _2}\left( t \right)} \right\|_{L^\infty} \leqslant C{\left( {1 + t} \right)^{ - \frac{3}
{2}}}{\left\| {\left[ {{\rho _{2,0}},{u_{2,0}},{\Theta _{2,0}}}
\right]} \right\|_{{L^1} \cap {{\dot H}^5}}}.
\end{split}\right.
\end{equation}
\end{coro}
%
\section{ Decay rates for system \eqref{2.2}}
\subsection{Decay rates for the energy functional.} In
this subsection, we will prove the decay rate \eqref{2.12} in
Proposition \ref{prop2.2} for the energy
$\left\|{U(t)}\right\|^2_s$. We begin with the Lemma as follows
which can be seen directly from the proof of Theorem \ref{thm2.1}.
\begin{lemma}\label{L4.1}
 Assume that $U=[\rho_\mu,~u_\mu,~\Theta_\mu,~E,~B]$ is the solution of the
 initial value
 problem \eqref{2.2}-\eqref{2.3} with $U_0$ $=[\rho_{\mu0}$, $u_{\mu0}$, $\Theta_{\mu0}$, $E_0,~B_0]$
 which satisfies \eqref{2.4}. If $\mathcal {E}_s(U_0)$
is small enough, then, for any $t \geq 0$
\begin{equation}\label{4.1}
\frac{d}{dt}\mathcal {E}_s(U(t))+\lambda \mathcal {D}_s(U(t))\leq 0.
\end{equation}
\end{lemma}
From Lemma \ref{L4.1}, we can check that
\begin{equation}\notag
\begin{split}
(1+t)^l\mathcal {E}_s(U(t))+& \gamma\int_0^t(1+y)^l\mathcal
{D}_s(U(y))dy\\
&\leq \mathcal {E}_{s}(U_0)+l\int_0^t(1+y)^{l-1}\mathcal
{E}_s(U(y)) dy\\
&\leq \mathcal {E}_{s}(U_0)+C l\int_0^t(1+y)^{l-1}\left(\left\|{B(y)
}\right\|^2+\left\|{(\rho_e+\rho_i)(y) }\right\|^2+\mathcal
{D}_{s+1}(U(y)) \right)dy,
\end{split}
\end{equation}
where we have used $\mathcal {E}_s(U(t))\leq \left\|{B(t)
}\right\|^2+\left\|{(\rho_e+\rho_i)(t) }\right\|^2+\mathcal
{D}_{s+1}(U(t))$. Using \eqref{4.1} again, we have
 $$\mathcal {E}_{s+2}(U(t))+\gamma\int_0^t \mathcal {D}_{s+2}(U(y))dy\leq \mathcal {E}_{s+2}(U_0)$$
and
\begin{equation}\notag
\begin{split}
(1+t)^{l-1}\mathcal {E}_{s+1}(U(t))+&
\gamma\int_0^t(1+y)^{l-1}\mathcal
{D}_{s+1}(U(y))dy\\
\leq \mathcal
{E}_{s+1}(U_0)&+C(l-1)\int_0^t(1+y)^{l-2}\left(\left\|{B(y)
}\right\|^2+\left\|{( \rho_e+\rho_i ) (y) }\right\|^2+\mathcal
{D}_{s+2}(U(y)) \right)dy.
\end{split}
\end{equation}
Then, by iterating the previous estimates, we obtain
\begin{equation}\label{4.2}
\begin{split}
(1+t)^l\mathcal {E}_s(U(t))+& \gamma\int_0^t(1+y)^l\mathcal
{D}_s(U(y))dy\\
&\leq C\mathcal
{E}_{s+2}(U_0)+C\int_0^t(1+y)^{l-1}\left(\left\|{B(y)
}\right\|^2+\left\|{( \rho_e+\rho_i ) (y) }\right\|^2\right)dy
\end{split}
\end{equation}
for $1<l<2.$

 Now, let us establish the estimate on the integral term on the right
hand side of \eqref{4.2}. Applying the estimate on $B$ in
\eqref{3.12} and the estimate on $\rho_2$ in \eqref{3.43} to
\eqref{3.5} and \eqref{3.6}, respectively,
 we have \begin{equation}\label{4.3}
 \begin{split}\left\| {B\left( t \right)}
\right\| \leqslant & C {\left( {1 + t} \right)^{ - \frac{3}
{4}}}{\left\| {\left[ {{u_{1,0}},{E_0},{B_0}} \right]}
\right\|_{{L^1}
\cap {{\dot H}^2}}}\\
& + C\int_0^t {{{\left( {1 + t - y} \right)}^{ - \frac{3}
{4}}}{{\left\| {\left[
{{g_{2e}}(y)-{g_{2i}}(y)},~{{g_{4e}}(y)-{g_{4i}}(y)} \right]}
\right\|}_{{L^1} \cap {{\dot H}^2}}}} dy, \end{split}\end{equation}
\begin{equation}\label{4.4}
 \begin{split}
  \left\| {\left( {{\rho _e} + {\rho _i}} \right)\left( t \right)}
   \right\| \leqslant & C\left\| {{\rho _2}\left( t \right)} \right\| \leqslant  C{\left( {1 + t} \right)^{ - \frac{3}
{4}}}{\left\| {\left[ {{\rho _{\mu 0}},{u_{\mu 0}},{\Theta _{\mu 0}}} \right]} \right\|_{{L^1} \cap {{\dot H}^2}}} \\
   &\quad + C\int_0^t {{{\left( {1 + t - y} \right)}^{ - \frac{3}
{4}}}{{\left\| {\left[ {{g_{1e}} + {g_{1i}},{g_{2e}} +
{g_{2i}},{g_{3e}} + {g_{3i}}} \right](y)} \right\|}_{{L^1}
 \cap {{\dot H}^2}}}dy.}
\end{split}\end{equation}
 It is directly to check that for any $0\leq y \leq t$,
 \[{\left\| {\left[
{{g_{2e}}(y)-{g_{2i}}(y)},~{{g_{4e}}(y)-{g_{4i}}(y)} \right]}
\right\|_{{L^1} \cap {{\dot H}^2}}}
 \leqslant C{\mathcal {E}_s}\left( {U(y)} \right) \leqslant
  C{\left( {1 + y} \right)^{ - \frac{3}
{2}}}{\mathcal {E} _{s,\infty }}\left( {U(t)} \right),\]
\[{\left\| {\left[ {{g_{1e}} + {g_{1i}},{g_{2e}} + {g_{2i}},{g_{3e}} + {g_{3i}}} \right](y)}
\right\|_{{L^1} \cap {{\dot H}^2}}}
 \leqslant C{\mathcal {E}_s}\left( {U(y)} \right) \leqslant
  C{\left( {1 + y} \right)^{ - \frac{3}
{2}}}{\mathcal {E} _{s,\infty }}\left( {U(t)} \right),\]
 where ${\mathcal {E} _{s,\infty }}\left(
{U(t)} \right): = \mathop {\sup }\limits_{0 \leqslant y \leqslant t}
{\left( {1 + y} \right)^{\frac{3} {2}}}{\mathcal {E}_s}\left( {U(y)}
\right).$ Plugging the two previous inequalities into \eqref{4.3}
and \eqref{4.4} respectively implies
\begin{equation}\label{4.5}
\left\| {B\left( t \right)} \right\| \leqslant C {\left( {1 + t}
\right)^{ - \frac{3} {4}}}\left( {{{\left\| {\left[
{{u_{\mu0}},{E_0},{B_0}} \right]} \right\|}_{{L^1} \cap {{\dot
H}^2}}} + {\mathcal {E} _{s,\infty }}\left( {U(t)} \right) }
\right)\end{equation}
and
\begin{equation}\label{4.6}
\left\| {\left( {{\rho _e} + {\rho _i}} \right)\left( t \right)}
   \right\| \leqslant C {\left( {1 + t}
\right)^{ - \frac{3} {4}}}\left( {{{\left\| {\left[
{{\rho_{\mu0}},{u_{\mu0}},{\Theta_{\mu0}}} \right]} \right\|}_{{L^1}
\cap {{\dot H}^2}}} + {\mathcal {E} _{s,\infty }}\left( {U(t)}
\right) } \right).\end{equation}

 Next, we search the uniform bound of
${\mathcal {E} _{s,\infty }}\left( {U(t)} \right)$ which will imply
the decay rates of ${\mathcal {E} _{s }}\left( {U(t)} \right)$ or
equivalent to $\|U(t)\|_s^2$. In fact, after choosing
$l=\frac{3}{2}+\varepsilon$ in \eqref{4.2} with $\varepsilon>0$
sufficiently small and utilizing \eqref{4.5} and \eqref{4.6}, one
obtain
\begin{equation}\notag\begin{split}
  {\left( {1 + t} \right)^{\frac{3}
{2} + \varepsilon }}{\mathcal {E} _s} & \left( {U(t)} \right) +
\gamma {\int_0^t {\left( {1 + y} \right)} ^{\frac{3}
{2} + \varepsilon }}{\mathcal {D}_s}\left( {U(y)} \right)dy  \\
 &  \leqslant C{\mathcal {E} _{s + 2}}\left( {{U_0}} \right) +
C{\left( {1 + t} \right)^\varepsilon }\left( {\left\| {\left[
{{\rho_{\mu0}},{u_{\mu0}},{\Theta_{\mu0}},{E_0},{B_0}} \right]}
\right\|_{_{{L^1} \cap {{\dot H}^2}}}^2 + {{\left[ {{\mathcal {E}
_{s,\infty }}\left( {U(t)} \right)} \right]}^2}} \right),
\end{split} \end{equation}
which implies
$${\left(
{1 + t} \right)^{\frac{3} {2}}}{\mathcal {E}_s}\left( {U(t)} \right)
\leqslant C\left( {{\mathcal {E}_{s + 2}}\left( {{U_0}} \right) +
\left\| {\left[
{{\rho_{\mu0}},{u_{\mu0}},{\Theta_{\mu0}},{E_0},{B_0}} \right]}
\right\|_{{L^1}}^2 + {{\left[ {{\mathcal {E} _{s,\infty }}\left(
{U(t)} \right)} \right]}^2}} \right),
$$
 and therefore,
$$
{\mathcal {E} _{s,\infty }}\left( {U(t)} \right) \leqslant C\left(
{{\omega _{s + 2}}
  {{\left( {{U_0}} \right)}^2} + {{\left[ {{\mathcal {E}_{s,\infty }}\left( {U(t)} \right)}
  \right]}^2}} \right), $$
since $\omega _{s + 2}\left( {{U_0}} \right)>0 $ is small enough, it
holds that ${\mathcal {E} _{s,\infty }}\left( {U(t)} \right)
\leqslant C{\omega _{s + 2}}{\left( {{U_0}} \right)^2} $ for any
$t\geq 0,$ which implies ${\left\| {U(t)} \right\|_s} \leqslant
C{\mathcal {E}_s}{\left( {U(t)} \right)^{\frac{1} {2}}} \leqslant
C{\omega _{s + 2}}\left( {{U_0}} \right){\left( {1 + t} \right)^{ -
\frac{3} {4}}}$, that is \eqref{2.12}.
 \subsection{Decay rate for high-order energy
functional.} In this subsection, we will look for the decay estimate
of the high-order energy $\|\nabla U(t)\|^2_{s-1}$, that is
\eqref{2.13} of Proposition \ref{prop2.2}. We begin with the
following Lemma.
\begin{lemma}\label{L4.2}
Assume $U=[\rho_\mu,~u_\mu,~\Theta_\mu,~E,~B]$ is the solution of
the initial value problem \eqref{2.2}-\eqref{2.3} with $U_0$
$=[\rho_{\mu0}$, $u_{\mu0}$, $\Theta_{\mu0}$, $E_0$, $B_0]$ which
satisfies \eqref{2.4} in the sense of Proposition \ref{prop2.1}. If
$\mathcal {E}_s(U_0)$ is small enough, then, there exist the
high-order energy function $\mathcal {E}_s^h(\cdot) $ and the
high-order dissipative rate $\mathcal {D}_s^h(\cdot) $ such that
\begin{equation}\label{4.7}
\frac{d}{dt}\mathcal {E}_s^h(U(t))+\gamma \mathcal {D}_s^h(U(t))\leq
0,
\end{equation}
holds for any $t\geq 0.$
\end{lemma}
\noindent \emph{Proof.}   The proof is very similar to the proof of
Theorem \ref{thm2.1}.  In fact, by letting $|\alpha|\geq 1$, then
corresponding to \eqref{2.16}, \eqref{2.19}, \eqref{2.23} and
\eqref{2.25},  it can also be checked that
\begin{equation}\notag
\frac{d} {{dt}}\left\|\nabla U \right\|_{s-1}^2 + \left\|
{\nabla\left[ {{u_e},{u_i},{\Theta _e},{\Theta _i}} \right]}
\right\|_{s-1}^2 \leqslant C{\left\| U \right\|_s} { \left\| {\nabla
\left[ {{\rho _e},{\rho _i},{u_e},{u_i},{\Theta _e},{\Theta _i}}
\right]} \right\|_{s - 1}^2},
\end{equation}
\begin{equation}\notag
\begin{split}
  \frac{d}
{{dt}}  \sum\limits_{1 \leq\left| \alpha  \right| \leqslant s - 1}
{\sum\limits_{\mu  = e,i} {\left\langle {{\partial ^\alpha }{u_\mu
},\nabla {\partial ^\alpha }{\rho _\mu }} \right\rangle } }
 &+ \gamma \left( {\left\| {\nabla^2 \left[ {{\rho _e},{\rho _i}} \right]} \right\|_{s - 2}^2
 + {{\left\| \nabla[{{\rho _e} - {\rho _i}}] \right\|}^2}} \right) \\
   & \leqslant C {\left( {\left\| {\nabla{u_\mu }} \right\|_{s-1}^2
   + \left\| {U} \right\|_s^2 {\left\|
    {\nabla [{\rho _\mu },{u_\mu },{\Theta _\mu }]} \right\|_{s - 1}^2 }
     } \right)},
\end{split}\end{equation}
\begin{equation}\notag
  \begin{split}
  \frac{d}
{{dt}}&\sum\limits_{1\leq\left| \alpha  \right| \leqslant s - 1}
{\left\langle {{\partial ^\alpha }\left( {{u_e}
 - {u_i}} \right),{\partial ^\alpha }E} \right\rangle }  + \gamma \left\| \nabla E \right\|_{s - 2}^2  \\
  & \leqslant C   \left( {\left\| {\nabla\left[ {{u_\mu },{\Theta _\mu }} \right]} \right\|_{s-1}^2
  + \left\|\nabla^2 \rho_\mu \right\|_{s-2}^2
  + \left\| {\nabla{u_\mu }}
   \right\|_{s-1}{{\left\| {\nabla^2 B} \right\|}_{s - 3}}} + \left\| U \right\|_s^2 {\left\|
   \nabla [{\rho _\mu },{u_\mu },{\Theta _\mu }] \right\|_{s - 1}^2 }
      \right),
        \end{split}
\end{equation}
and
\begin{equation}\notag
\begin{split}
\frac{d} {{dt}}\sum\limits_{1\leq\left| \alpha  \right| \leqslant s
- 2} &  {  \left\langle {{\partial ^\alpha }E, - \nabla  \times
{\partial ^\alpha }B} \right\rangle  +  \gamma \left\| {\nabla^2 B}
\right\|_{s - 3}^2}\\
& \leqslant  C  {(\left\| { {\nabla E} } \right\|_{s - 2}^2 +
\left\| {\nabla {{u_\mu }} } \right\|_{s - 1}^2 } +\left\| {\nabla
[{\rho _\mu },{u _\mu }]} \right\|_{s- 1}^2\left\| {{U }}
\right\|_s^2 ).
\end{split}\end{equation}

 Now, similarly done as that in
\emph{Step 5} of Theorem \ref{thm2.1}. Let us define the high-order
energy functional as
\begin{equation}\label{4.8}
  \begin{split}
 \mathcal{E}_s(U(t))=& \left\|\nabla U \right\|_{s-1}^2 +\mathcal {K}_1 \sum\limits_{1\leq \left| \alpha  \right|
  \leqslant s - 1}
  {\sum\limits_{\mu  = e,i}
   {\left\langle {{\partial ^\alpha }{u_\mu },\nabla {\partial ^\alpha }{\rho _\mu }} \right\rangle } }   \\
    &
+\mathcal {K}_2 \sum\limits_{1\leq \left| \alpha  \right| \leqslant
s - 1} {\left\langle {{\partial ^\alpha }\left( {{u_e}
    - {u_i}} \right),{\partial ^\alpha }E} \right\rangle }
+\mathcal {K}_3 \sum\limits_{1\leq \left| \alpha  \right| \leqslant
s - 2}
  {\left\langle {{\partial ^\alpha }E, - \nabla  \times {\partial ^\alpha }B} \right\rangle } ,
  \end{split}
\end{equation}
Similarly, one can take $0<\mathcal {K}_3\ll\mathcal
{K}_2\ll\mathcal {K}_1\ll 1$ be sufficiently small with $\mathcal
{K}_2^{\frac{3}{2}}\ll\mathcal {K}_3$, such that
$\mathcal{E}_s^h(U(t))\sim \|\nabla U(t)\|_{s-1}^2 $, that is
$\mathcal{E}_s^h(\cdot)$ is really a high-order energy functional
which satisfies \eqref{2.6},  and moreover, the sum of the four
previously  estimates with coefficients corresponding to \eqref{4.8}
gives \eqref{4.7}. Now, we complete the proof of Lemma \ref{L4.2}.
\hfill $\Box$

Based on Lemma \ref{L4.2}, one can check that
\begin{equation}\notag
\frac{d}{dt}\mathcal {E}_s^h(U(t))+\gamma \mathcal {E}_s^h(U(t))\leq
C\left(\|\nabla
B\|^2+\|\nabla^s[E,B]\|^2+\|\nabla(\rho_e+\rho_i)\|^2\right),
\end{equation}
which implies
\begin{equation}\label{4.9}
\begin{split}
\mathcal {E}_s^h(U(t))\leq & e^{-\gamma t}\mathcal {E}_s^h(U_0)
\\&
+ C\int_0^t{e^{-\gamma (t-y)}\left(\|\nabla B(y)\|^2+\|\nabla^s
[E,~B](y) \|^2 + \|\nabla (\rho_e+\rho_i)(y) \|^2 \right)}dy.
\end{split}
\end{equation}
 Now, let us estimate the time integral term on the right hand side of the
previous inequality. Noting that the equations of $E$ and $B$ in
bipolar non-isentropic Euler-Maxwell system are the same as that in
bipolar isentropic Euler-Maxwell system, similarly as that in
\cite{Duan11b}, we obtain
\begin{lemma}\label{L4.3}
Assume $U=[\rho_\mu,~u_\mu,~\Theta_\mu,~E,~B]$ is the solution of
the initial value problem \eqref{2.2}-\eqref{2.3} with $U_0$
$=[\rho_{\mu0}$, $u_{\mu0}$, $\Theta_{\mu0}$, $E_0$, $B_0]$ which
satisfies \eqref{2.4} in the sense of Proposition \ref{prop2.1}. If
$\omega_{s+6}(U_0)$ is small enough, then, for any $t\geq 0$
\begin{equation}\label{4.10}
{\left\| {\nabla B(t)} \right\|^2} + {\left\| {{\nabla ^s}\left[
{E(t),B(t)} \right]} \right\|^2}+\|\nabla (\rho_e+\rho_i)(t) \|^2
\leqslant C{\omega _{s + 6}}{({U_0})^2}{(1 + t)^{ - \frac{5} {2}}}.
\end{equation}
\end{lemma}
\noindent \emph{Proof.} Utilize the estimate \eqref{3.14} to
\eqref{3.5} of the solution $U_1(t)$ so that
\begin{equation}\notag
\begin{split}
  \left\| {\nabla B\left( t \right)} \right\| & \leqslant C{\left( {1 + t} \right)^{ - \frac{5}
{4}}}{\left\| {\left[ {{u_{\mu0}},{E_0},{B_0}} \right]}
\right\|_{{L^1} \cap {{\dot H}^4}}}\\
&\quad + C\int_0^t{\left( {1 +t- y} \right)^{ - \frac{5}
{4}}}{\left\| {\left[ {{g_{2e}(y)-g_{2i}(y)},{g_{4e}(y)-g_{4i}(y)}} \right]} \right\|_{{L^1} \cap {{\dot H}^4}}}dy \\
 & \leqslant C{\left( {1 + t} \right)^{ - \frac{5}
{4}}}{\left\| {\left[ {{u_{\mu0}},{E_0},{B_0}} \right]}
\right\|_{{L^1} \cap {{\dot H}^4}}}+ C\int_0^t{\left( {1 +t- y}
\right)^{ - \frac{5}
{4}}}{\left\| {U(y)} \right\|^2_{\max\{5,s\}  }}dy\\
 & \leqslant C{\left( {1 + t} \right)^{ - \frac{5}
{4}}}{\left\| {\left[ {{u_{\mu0}},{E_0},{B_0}} \right]}
\right\|_{{L^1} \cap {{\dot H}^4}}}+ C\int_0^t{\left( {1 +t- y}
\right)^{ - \frac{5} {4}}}\omega_{s+6}(U_0)^2   \left( {1 + y}
\right)^{ - \frac{3}
{2}} dy\\
 & \leqslant C\omega_{s+6}(U_0)   \left( {1 + t} \right)^{ - \frac{5}
{4}}
\end{split}
\end{equation}
and
\begin{equation}\notag
\begin{split}
  & \left\| {{\nabla ^s}\left[ {E\left( t \right),B\left( t \right)} \right]}
  \right\|\\
  & \leqslant
  C{\left( {1 + t} \right)^{ - \frac{5}
{4}}}{\left\| {\left[ {{u_{\mu0}},\Theta_{\mu0},{E_0},{B_0}}
\right]} \right\|_{{L^2} \cap {{\dot H}^{s + 3}}}}\\&
\quad+C\int_0^t{\left( {1 + t-y} \right)^{ - \frac{5} {4}}}{\left\|
{\left[
{{g_{2e}(y)-g_{2i}(y)},g_{3e}(y)-g_{3i}(y),{g_{4e}(y)-g_{4i}(y)}}
\right]} \right\|_{{L^2} \cap {{\dot
H}^{s + 3}}}}dy\\
  & \leqslant
  C{\left( {1 + t} \right)^{ - \frac{5}
{4}}}{\left\| {\left[ {{u_{\mu0}},\Theta_{\mu0},{E_0},{B_0}}
\right]} \right\|_{{L^2} \cap {{\dot H}^{s + 3}}}}+C\int_0^t{\left(
{1 + t-y} \right)^{ - \frac{5} {4}}}{\left\| { {U(y)} } \right\|^2_{
{{s +
4}}}}dy\\
  & \leqslant
  C{\left( {1 + t} \right)^{ - \frac{5}
{4}}}{\left\| {\left[ {{u_{\mu0}},\Theta_{\mu0},{E_0},{B_0}}
\right]} \right\|_{{L^2} \cap {{\dot H}^{s + 3}}}}+C\int_0^t{\left(
{1 + t-y} \right)^{ - \frac{5} {4}}}\omega_{s+6}(U_0)^2   \left( {1
+ y} \right)^{ - \frac{3}
{2}} dy\\
 & \leqslant C\omega_{s+6}(U_0)   \left( {1 + t} \right)^{ - \frac{5}
{4}}.
\end{split}
\end{equation}

 Similarly, utilizing the estimate on $\rho_2$ in
\eqref{3.44} to \eqref{3.6} of the solution $U_2(t)$, we obtain
\begin{equation}\notag
\begin{split}
  & \left\| {{\nabla }\left( {\rho_e+\rho_i} \right)\left( t \right)}
  \right\|\\
  & \leqslant
  C{\left( {1 + t} \right)^{ - \frac{5}
{4}}}{\left\| {\left[ {\rho_{\mu0},{u_{\mu0}},\Theta_{\mu0}}
\right]} \right\|_{{L^1} \cap {{\dot H}^{4}}}}\\&
\quad+C\int_0^t{\left( {1 + t-y} \right)^{ - \frac{5} {4}}}{\left\|
{\left[
{g_{1e}(y)+g_{1i}(y)},{{g_{2e}(y)+g_{2i}(y)},g_{3e}(y)+g_{3i}(y)}
\right]} \right\|_{{L^1} \cap {{\dot
H}^{4}}}}dy\\
  & \leqslant
  C{\left( {1 + t} \right)^{ - \frac{5}
{4}}}{\left\| {\left[ {\rho_{\mu0},{u_{\mu0}},\Theta_{\mu0}}
\right]} \right\|_{{L^1} \cap {{\dot H}^{4}}}}+C\int_0^t{\left(
{1 + t-y} \right)^{ - \frac{5} {4}}}{\left\| { {U(y)} } \right\|^2_{\max\{5,s\}  }}dy\\
  & \leqslant
  C{\left( {1 + t} \right)^{ - \frac{5}
{4}}}{\left\| {\left[ {\rho_{\mu0},{u_{\mu0}},\Theta_{\mu0}}
\right]} \right\|_{{L^2} \cap {{\dot H}^{s + 3}}}}+C\int_0^t{\left(
{1 + t-y} \right)^{ - \frac{5} {4}}}\omega_{s+6}(U_0)^2   \left( {1
+ y} \right)^{ - \frac{3}
{2}} dy\\
 & \leqslant C\omega_{s+6}(U_0)   \left( {1 + t} \right)^{ - \frac{5}
{4}}.
\end{split}
\end{equation}
Where we have used \eqref{2.12} and the smallness of
$\omega_{s+6}(U_0)  $. Now, we complete the proof of Lemma
\ref{L4.3}. \hfill $\Box$\\

 Then, after plugging \eqref{4.10} into \eqref{4.9},
we have
\begin{equation}\notag
\mathcal {E}_s^h(U(t))\leq e^{-\gamma t}\mathcal {E}_s^h(U_0)+ C
\omega_{s+6}(U_0)^2(1+t)^{-\frac{5}{2}}.
\end{equation}
Since $\mathcal {E}^h_s(U(t))\sim \|\nabla U(t) \|^2_{s-1}$ holds
true for any $t\geq 0$, \eqref{2.13} follows. Now, we finish the
proof of Proposition \ref{2.2}.
 \subsection{Decay rate in $L^q$ .}  In this subsection, we are to look for the
decay rates of solutions $U$ $=[\rho_\mu$, $u_\mu$, $\Theta_\mu$,
$E,~B]$ in $L^q$ $_{(2\leq q\leq +\infty)}$ of the initial value
problem \eqref{2.2}-\eqref{2.3} by proving the second part of
Theorem \ref{thm1.1}. Throughout this subsection, we usually suppose
that $\omega_{13}(U_0)$ is small enough. Firstly, for $s\geq 4$,
Proposition \ref{prop2.2} shows that if $\omega_{s+2}(U_0)$ is small
enough,
\begin{equation}\label{4.11}
\|U(t)\|_s\leq C \omega_{s+2}(U_0)(1+t)^{-\frac{3}{4}},
\end{equation}
and if $\omega_{s+6}(U_0)$ is small enough,
\begin{equation}\label{4.12}
\|\nabla U(t)\|_{s-1}\leq C \omega_{s+6}(U_0)(1+t)^{-\frac{5}{4}}.
\end{equation}
Now, let us establish the estimates on $B$, $[u_e-u_i,~E]$,
$u_e+u_i$, $[\rho_e-\rho_i,~\Theta_e-\Theta_i]$ and
$[\rho_e+\rho_i,~\Theta_e+\Theta_i]$
as follows.\\
 \noindent \emph{Estimate on $\|B\|_{L^q} $.} For $L^2$ rate, it is
 directly from \eqref{4.11} to get
 $$\|B(t)\|\leq C \omega_{6} (U_0)(1+t)^{-\frac{3}{4}}.$$
 For $L^\infty$ rate, by applying $L^\infty$ estimate on $B$
 of \eqref{3.13} to \eqref{3.5}, we obtain
\begin{equation}\notag
\begin{split}
{\left\| {B(t)} \right\|_{{L^\infty }}} \leqslant & C{(1 + t)^{ -
\frac{3} {2}}}{\left\| {[{u_{\mu0}},{E_0},{B_0}]} \right\|_{{L^1}
\cap {{\dot H}^5}}} \\
& + C\int_0^t {{{(1 + t - y)}^{ - \frac{3} {2}}}} {\left\|
{[{g_{2e}}-{g_{2i}},{g_{4e}}-{g_{4i}}](y)} \right\|_{{L^1} \cap
{{\dot H}^5}}}dy.
\end{split}
\end{equation}
 Because of \eqref{4.11},
\begin{equation}\notag
\begin{split}
{\left\| {[{g_{2e}}-{g_{2i}},{g_{4e}}-{g_{4i}}](t)} \right\|_{{L^1}
\cap {{\dot H}^5}}} \leqslant C\left\|
 {U(t)} \right\|_6^2 \leqslant C{\omega _8}{({U_0})^2}{(1 + t)^{ - \frac{3}
{2}}},\end{split}
\end{equation}
we have
 $${\left\| {B(t)} \right\|_{{L^\infty }}} \leqslant C\omega_8 ({U_0}) {(1 + t)^{ -
 \frac{3}
{2}}}.$$ Therefore, by $L^2-L^\infty$ interpolation
\begin{equation}\label{4.13}
\|B(t)\|_{L^q}\leq C \omega_8 ({U_0}) {(1 + t)^{ - \frac{3}
{2}+\frac{3}{2q}}},
\end{equation}
for $2\leq q \leq \infty.$\\
\noindent\emph{Estimate on $\|[u_e-u_i,E]\|_{L^q}$.} For $L^2$ rate,
applying the $L^2$ estimate on $u_e-u_i$ and $E$ in \eqref{3.12} to
\eqref{3.5}, one has
\begin{equation}\notag
\begin{split}
  \left\| {(u_e-u_i)\left( t \right)} \right\| \leqslant & C{\left( {1 + t} \right)^{ - \frac{5}
{4}}}\left( {\left\| {\left[ {{\rho _{\mu0}},{\Theta _{\mu0}}}
\right]} \right\| + {{\left\|
{\left[ {{u_{\mu0}},{E_0},{B_0}} \right]} \right\|}_{{L^1} \cap {{\dot H}^2}}}} \right)  \\
    & + C\int_0^t {{\left( {1 + t - y} \right)}^{ - \frac{5}
{4}}}  \left\| {\left[ {{g_{1e}-g_{1i}},{g_{3e}-g_{3i}}} \right](y)}
\right\|dy
 \\ &  + C\int_0^t {{\left( {1 + t - y} \right)}^{ - \frac{5}
{4}}}  {{\left\| {\left[ {{g_{2e}-g_{2i}},{g_{4e}-g_{4i}}}
\right](y)} \right\|}_{{L^1} \cap {{\dot H}^2}}}
  dy
\end{split}
\end{equation}
and
\begin{equation}\notag
\begin{split}
\left\| {E\left( t \right)} \right\| \leqslant & C{\left( {1 + t}
\right)^{ - \frac{5} {4}}}{\left\| {\left[
{{u_{\mu0}},\Theta_{\mu0},{E_0},{B_0}} \right]} \right\|_{{L^1} \cap
{{\dot
H}^3}}}\\
& + C\int_0^t {{{\left( {1 + t - y} \right)}^{ - \frac{5}
{4}}}{{\left\| {\left[
{{g_{2e}-g_{2i}},{g_{3e}-g_{3i}},{g_{4e}-g_{4i}}} \right](y)}
\right\|}_{{L^1} \cap {{\dot H}^3}}}dy.}
\end{split}
\end{equation}
 Since by \eqref{4.11},
\begin{equation}\notag
\begin{split}
\left\| {\left[ {{g_{1e}-g_{1i}},  {g_{3e}-g_{3i}}} \right](t)}
\right\| + {\left\| {\left[
{{g_{2e}-g_{2i}},{g_{3e}-g_{3i}},{g_{4e}-g_{4i}}}
\right](t)} \right\|_{{L^1} \cap {{\dot H}^3}}}\\
 \leqslant C\left\| {U(t)} \right\|_4^2  \leqslant C{\omega
_6}{({U_0})^2}{(1 + t)^{ - \frac{3} {2}}},\end{split}
\end{equation}
 which implies that
\begin{equation}\label{4.14}
 \left\|  [u_e-u_i,~E]( t )\right\|\leqslant C{\omega _6}({U_0}){\left( {1 + t} \right)^{ -
\frac{5} {4}}}.
\end{equation}
 For $L^\infty$ rate, utilize the $L^\infty$ estimates on $u_e-u_i$ and $E$ in \eqref{3.13} to \eqref{3.5},
 we have
\begin{equation}\notag
\begin{split}
  {\left\| {(u_e-u_i)\left( t \right)} \right\|_{{L^\infty }}} \leqslant& C{\left( {1 + t} \right)^{ - 2}}
  \left( {{{\left\| {\left[ {{\rho _{\mu0}},{\Theta _{\mu0}}} \right]} \right\|}_{{L^1} \cap {{\dot H}^2}}}
   + {{\left\| {\left[ {{u_{\mu0}},{E_0},{B_0}} \right]} \right\|}_{{L^1} \cap {{\dot H}^5}}}} \right)  \\
    & + C\int_0^t {{\left( {1 + t - y} \right)}^{ -2}}  \left\| {\left[ {{g_{1e}-g_{1i}},{g_{3e}-g_{3i}}}
\right](y)} \right\|_{{L^1} \cap {{\dot H}^2}}dy
 \\ &  + C\int_0^t {{\left( {1 + t - y} \right)}^{ - 2}}  {{\left\| {\left[ {{g_{2e}-g_{2i}},{g_{4e}-g_{4i}}}
\right](y)} \right\|}_{{L^1} \cap {{\dot H}^5}}}
  dy
\end{split} \end{equation}
and
\begin{equation}\notag
\begin{split}{\left\| {E\left( t \right)} \right\|_{{L^\infty }}}
\leqslant & C{\left( {1 + t} \right)^{ - 2}}{\left\| {\left[
{{u_{\mu0}},{\Theta_{\mu0}},{E_0},{B_0}} \right]} \right\|_{{L^1}
\cap {{\dot
H}^6}}}\\
& + C\int_0^t {{{\left( {1 + t - y} \right)}^{ - 2}}{{\left\|
{\left[ {{g_{2e}-g_{2i}},{g_{3e}-g_{3i}},{g_{4e}-g_{4i}}}
\right](y)} \right\|}_{{L^1} \cap {{\dot H}^6}}}dy}.
\end{split} \end{equation}
 Since
\begin{equation}\notag
\begin{split} \|[ {g_{1e}-g_{1i}},{g_{2e}-g_{2i}},&{g_{3e}-g_{3i}},{g_{4e}-g_{4i}}
](t)\|_{L^1} \\& \leq C \|U(t)\|(\|(u_e-u_i)(t)\|+\|\ U(t)\|+\|\nabla U(t)\|)\\
&\leq \omega_{10}(U_0)^2(1+t)^{-\frac{3}{2}},
\end{split} \end{equation}
and
\begin{equation}\notag
\begin{split}
  \|[ {g_{1e}-g_{1i}},{g_{2e}-g_{2i}},&{g_{3e}-g_{3i}},{g_{4e}-g_{4i}}
](t)\|_{{{\dot H}^5} \cap {{\dot H}^6}}
    \leqslant C\left\| {\nabla U(t)} \right\|_6^2 \leqslant {\omega _{13}}{({U_0})^2}{(1 + t)^{ - \frac{5}
{2}}},
\end{split} \end{equation}
 then, one has
 $$\|[u_e(t)-u_i(t),~E(t)]\|_{L^\infty}\leq C \omega_{13}(U_0)^2(1+t)^{-\frac{3}{2}}.$$
 Therefore, by $L^2-L^\infty$ interpolation
\begin{equation}\label{4.15}
\|[u_e(t)-u_i(t),~E(t)]\|_{L^q}\leq C \omega_{13} ({U_0}) {(1 + t)^{
- \frac{3}{2}+\frac{1}{2q}}},
\end{equation}
for $2\leq q \leq \infty.$

\noindent\emph{Estimate on $\|u_e+u_i\|_{L^q}$.} For $L^2$ rate,
utilizing the $L^2$ estimates on $u_e+u_i$ in \eqref{3.43} to
\eqref{3.6}, we have
\begin{equation}\notag
\begin{split}
  \left\| {(u_e+u_i)\left( t \right)} \right\| \leqslant & C{\left( {1 + t} \right)^{ - \frac{5}
{4}}} {\left\| {\left[ {{\rho _{\mu0}},{u _{\mu0}},{\Theta
_{\mu0}}} \right]} \right\|_{{L^1} \cap {{\dot H}^3}} }   \\
    & + C\int_0^t {{\left( {1 + t - y} \right)}^{ - \frac{5}
{4}}}  \left\| {\left[
{{g_{1e}+g_{1i}},g_{2e}+g_{2i},{g_{3e}+g_{3i}}} \right](y)}
\right\|_{{L^1} \cap {{\dot H}^3}}dy.
\end{split}
\end{equation}
 Since by \eqref{4.11},
\begin{equation}\notag
\begin{split} \left\| {\left[
{{g_{1e}+g_{1i}},g_{2e}+g_{2i},{g_{3e}+g_{3i}}} \right](t)}
\right\|_{{L^1} \cap {{\dot H}^3}}
 \leq C \|U(t)\|^2_4 \leq \omega_{6}(U_0)^2(1+t)^{-\frac{3}{2}},
\end{split} \end{equation}
it follows that
\begin{equation}\notag
 \left\|  (u_e+u_i)( t )\right\|\leqslant C{\omega _6}({U_0}){\left( {1 + t} \right)^{ -
\frac{5} {4}}}.
\end{equation}
For $L^\infty$ rate, utiliz the $L^\infty$ estimates on $u_e+u_i$ in
\eqref{3.45} to \eqref{3.6}, we have
\begin{equation}\notag
\begin{split}
  {\left\| {(u_e+u_i)\left( t \right)} \right\|_{{L^\infty }}} \leqslant& C{\left( {1 + t} \right)^{ - 2}}
  {\left\| {\left[ {{\rho _{\mu0}},{u _{\mu0}},{\Theta
_{\mu0}}} \right]} \right\|_{{L^1} \cap {{\dot H}^6}} } \\
    & + C\int_0^t {{\left( {1 + t - y} \right)}^{ -2}}  \left\| {\left[
{{g_{1e}+g_{1i}},g_{2e}+g_{2i},{g_{3e}+g_{3i}}} \right](y)}
\right\|_{{L^1} \cap {{\dot H}^6}}dy
\end{split} \end{equation}
Since by \eqref{4.11},
\begin{equation}\notag
\begin{split} \|[ {g_{1e}+g_{1i}},{g_{2e}+g_{2i}},&{g_{3e}+g_{3i}}
](t)\|_{{L^1} \cap {{\dot H}^6}} \leq C \|U(t)\|^2_7\leq
\omega_{9}(U_0)^2(1+t)^{-\frac{3}{2}},
\end{split} \end{equation}
 it follows that
 $$\|u_e(t)+u_i(t)\|_{L^\infty}\leq C \omega_{9}(U_0) (1+t)^{-\frac{3}{2}}.$$
 Therefore, by $L^2-L^\infty$ interpolation
\begin{equation}\label{4.16}
\|u_e(t)+u_i(t)\|_{L^q}\leq C \omega_{9} ({U_0}) {(1 + t)^{ -
\frac{3}{2}+\frac{1}{2q}}},
\end{equation}
for $2\leq q \leq \infty.$

Then from \eqref{4.15} and \eqref{4.16} we have
\begin{equation}\label{4.17}
\|u_\mu (t)\|_{L^q}\leq C \omega_{13} ({U_0}) {(1 + t)^{ -
\frac{3}{2}+\frac{1}{2q}}},
\end{equation}
for $2\leq q \leq \infty.$
%

\noindent\emph{Estimate on
$\|[\rho_e-\rho_i,\Theta_e-\Theta_i]\|_{L^q}$ and
$\|[\rho_e+\rho_i,\Theta_e+\Theta_i]\|_{L^q}$.} For $L^2$ rate,
utilizing the $L^2$ estimates on $\rho_e-\rho_i$ and
$\Theta_e-\Theta_i$ in \eqref{3.12} to \eqref{3.5}, we have
\begin{equation}\label{4.18}
\begin{split}
  & \left\| {\left[ {\rho_e-\rho_i}, {\Theta_e-\Theta_i}\right]\left( t \right)}
  \right\| \\
 &\quad  \leqslant
  C{e^{ - \frac{t}
{2}}}{\left\| {\left[ {\rho_{\mu0},{u_{\mu0}},\Theta_{\mu0}}
\right]} \right\|}   +C\int_0^t{e^{ - \frac{{t - y}} {2}}}{\left\|
{\left[ {g_{1e} -g_{1i} },{{g_{2e} -g_{2i} },g_{3e} -g_{3i}}
\right](y)} \right\|}dy.
\end{split}
\end{equation}
Because of
\begin{equation}\notag
\begin{split}
& \| {\left[ {g_{1e} -g_{1i} },{{g_{2e} -g_{2i} },g_{3e} -g_{3i}}
\right](t)} \|\\
& \leqslant  C\left(\left\| {\nabla U(t)}
  \right\|_1^2 + \left\| {(u_e+u_i)(t)} \right\|
  {\left\| { {B(t)} } \right\|_{{L^\infty }}}\right)
   \leqslant C{\omega _{10}}{({U_0})^2}{\left( {1 + t} \right)^{ - \frac{5}
{2}}},
 \end{split}
\end{equation}
where \eqref{4.12}, \eqref{4.13} and \eqref{4.16} were used. Then
\eqref{4.18} yields the decay estimate
\begin{equation}\label{4.19}
\begin{split}\left\| {\left[ {\rho_e-\rho_i}, {\Theta_e-\Theta_i}\right]\left( t \right)}
  \right\|\leqslant C
{\omega _{10}}{({U_0})}{\left( {1 + t} \right)^{ - \frac{5}
{2}}}.\end{split}
\end{equation}
Similarly for $\|[\rho_e-\rho_i,\Theta_e-\Theta_i]\|$,  by utilizing
the ${L^2}$ estimate on $[\rho_e+\rho_i,\Theta_e+\Theta_i]$ in
\eqref{3.43} to \eqref{3.6}, we obtain the decay estimate
\begin{equation}\label{4.20}
\begin{split}\left\| {\left[ {\rho_e+\rho_i}, {\Theta_e+\Theta_i}\right]\left( t \right)}
  \right\|\leqslant C
{\omega _{6}}{({U_0})}{\left( {1 + t} \right)^{ - \frac{3}
{4}}}.\end{split}
\end{equation}
 Combining \eqref{4.19} and \eqref{4.20}, we obtain
\begin{equation}\label{4.21}
\begin{split}\left\| {\left[ {\rho_\mu}, {\Theta_\mu}\right]\left( t \right)}
  \right\|\leqslant C
{\omega _{10}}{({U_0})}{\left( {1 + t} \right)^{ - \frac{3}
{4}}}.\end{split}
\end{equation}
 For $L^{\infty}$
rate, by utilizing the $L^{\infty}$ estimate on
$[\rho_e-\rho_i,\Theta_e-\Theta_i]$ in \eqref{3.13} to \eqref{3.5},
we have the decay estimate
\begin{equation}\label{4.22}
\begin{split}
   \left\| {\left[ {\rho_e-\rho_i}, {\Theta_e-\Theta_i}\right]\left( t \right)}
  \right\|_{L^\infty}  \leqslant &
  C{e^{ - \frac{t}
{2}}}{\left\| {\left[ {\rho_{\mu0},{u_{\mu0}},\Theta_{\mu0}}
\right]} \right\|_{L^2\cap\dot{H}^2}} \\
 &   +C\int_0^t{e^{ - \frac{{t -
y}} {2}}}{\left\| {\left[ {g_{1e} -g_{1i} },{{g_{2e} -g_{2i}
},g_{3e} -g_{3i}} \right](y)} \right\|_{L^2\cap\dot{H}^2}}dy.
\end{split}
\end{equation}
 Notice that one can check
\begin{equation}\label{4.23}
\begin{split}
  &{{\left\| {\left[ {g_{1e} -g_{1i} },{{g_{2e} -g_{2i}
},g_{3e} -g_{3i}} \right](t)} \right\|}_{{L^2}\cap {{\dot H}^2}}}\\
   &\quad\quad \leqslant C{\left\|
  {\nabla U(t)} \right\|_4}\left({\left\| {\left[ {\rho_\mu (t),\Theta_\mu (t)} \right]} \right\|}+
  \left\| {{u_\mu (t)} }
  \right\|+
   {{{\left\| {\left[ {u_\mu (t),B(t)} \right]} \right\|}_{{L^\infty }}} } \right)  \\
  &\quad\quad \leq C{\omega _{13}}{({U_0})^2}{(1 +
t)^{ - 2}},\end{split}
\end{equation}
where we have used \eqref{4.12}, \eqref{4.13}, \eqref{4.17} and
\eqref{4.21}. Which implies from \eqref{4.22} that
\[\left\| {\left[ {\rho_e-\rho_i}, {\Theta_e-\Theta_i}\right]\left( t \right)}
  \right\|_{L^\infty}  \leqslant  C{\omega _{13}}{({U_0})}(1 +
t)^{ - 2}.\]
 Therefore, by $L^2-L^\infty$ interpolation
\begin{equation}\label{4.24}
\|\left[ {\rho_e-\rho_i}, {\Theta_e-\Theta_i}\right]\|_{L^q}\leq C
\omega_{13} ({U_0}) {(1 + t)^{ -2- \frac{1}{q}}},
\end{equation}
for $2\leq q \leq \infty.$

For $\|\left[ {\rho_e+\rho_i},
{\Theta_e+\Theta_i}\right]\|_{L^\infty}$, by utilizing the
$L^\infty$ estimate on $ \left[ {\rho_e+\rho_i},
{\Theta_e+\Theta_i}\right] $ in \eqref{3.45} to \eqref{3.6}, we have
the decay estimate
\begin{equation}\label{4.25}
\begin{split}\left\| {\left[ {\rho_e+\rho_i}, {\Theta_e+\Theta_i}\right]\left( t \right)}
  \right\|_{L^\infty} \leqslant C
{\omega _{8}}{({U_0})}{\left( {1 + t} \right)^{ - \frac{3}
{2}}}.\end{split}
\end{equation}
Then from \eqref{4.20} and \eqref{4.25} we have
\begin{equation}\label{4.26}
\begin{split}\left\| {\left[ {\rho_e+\rho_i}, {\Theta_e+\Theta_i}\right]\left( t \right)}
  \right\|_{L^q} \leqslant C
{\omega _{8}}{({U_0})}{\left( {1 + t} \right)^{ - \frac{3}
{2}+\frac{3}{2q}}}.\end{split}
\end{equation}
 Thus, \eqref{4.24},  \eqref{4.26},  \eqref{4.15}-\eqref{4.16} and
 \eqref{4.13} give  \eqref{1.4},  \eqref{1.5},  \eqref{1.6} and  \eqref{1.7}, respectively.
 Now, we complete the proof of Theorem \ref{thm1.1}.\hfill $\Box$\\
\noindent {\bf Acknowledgments}
 \vspace{2mm}
This work is supported by  the NSFC (Grant no. 11071009), BSF (Grant
no. 1082001), the fund of Beijing education committee of China, and the Foundation
Project of Doctor Graduate Student Innovation of Beijing University
of Technology of China.


\end{document}